\documentclass[12pt]{article}
\usepackage[dvips]{epsfig}
\usepackage{graphicx}
\usepackage{amsthm,amsmath,amssymb}
\usepackage{mathtools}
\usepackage{latexsym}
\usepackage{graphics}
\usepackage{url}
\usepackage{placeins}
\usepackage{color}
\usepackage{bm}
\usepackage{comment}
\usepackage{statex}
\usepackage{pgfplots}

\usepackage[
            bookmarksnumbered=true,
            bookmarksopen=true,
            colorlinks=true,
            citecolor=green,
            linkcolor=green,
            anchorcolor=red,
            urlcolor=blue
            ]{hyperref}

\definecolor{blue}{rgb}{0,0,0.9}
\definecolor{red}{rgb}{0.9,0,0}
\definecolor{green}{rgb}{0,0.50,0.10}
\definecolor{violet}{rgb}{0.5804,0.0000,0.8275}

\def\@themcountersep{}

\makeatletter
\newcommand{\labeltext}[2]{%
  \@bsphack
  \csname phantomsection\endcsname 
  \def\@currentlabel{#1}{\label{#2}}%
  \@esphack
}
\makeatother

\newtheorem{THEO}{Theorem}[section]
\newtheorem{ALGo}[THEO]{Algorithm}
\newtheorem{ASSUMPT}[THEO]{Assumption}
\newtheorem{CONJ}[THEO]{Conjecture}
\newtheorem{COND}[THEO]{Condition}
\newtheorem{CORO}[THEO]{Corollary}
\newtheorem{DEFI}[THEO]{Definition}
\newtheorem{EXAMP}[THEO]{Example}
\newtheorem{INSTANCE}[THEO]{Instance}
\newtheorem{FACT}[THEO]{Fact}
\newtheorem{HYPO}[THEO]{Hypothesis}
\newtheorem{LEMM}[THEO]{Lemma}
\newtheorem{PROB}[THEO]{Problem}
\newtheorem{PROP}[THEO]{Proposition}
\newtheorem{REMA}[THEO]{Remark}
\newcommand{\theo}{\begin{THEO}}
\newcommand{\algo}{\begin{ALGo} \rm}
\newcommand{\assumpt}{\begin{ASSUMPT} \rm}
\newcommand{\cond}{\begin{COND}}
\newcommand{\conj}{\begin{CONJ}}
\newcommand{\coro}{\begin{CORO}}
\newcommand{\defi}{\begin{DEFI} \rm}
\newcommand{\examp}{\begin{EXAMP} \rm}
\newcommand{\instan}{\begin{INSTANCE} \rm}
\newcommand{\fact}{\begin{FACT}}
\newcommand{\hypo}{\begin{HYPO} \rm}
\newcommand{\lemm}{\begin{LEMM}}
\newcommand{\prob}{\begin{PROB} \rm}
\newcommand{\prop}{\begin{PROP}}
\newcommand{\rema}{\begin{REMA} \rm}

\newcommand{\etheo}{\end{THEO}}
\newcommand{\ealgo}{\end{ALGo}}
\newcommand{\eassumpt}{\end{ASSUMPT}}
\newcommand{\econd}{\end{COND}}
\newcommand{\econj}{\end{CONJ}}
\newcommand{\ecoro}{\end{CORO}}
\newcommand{\edefi}{\end{DEFI}}
\newcommand{\eexamp}{\end{EXAMP}}
\newcommand{\einstan}{\end{INSTANCE}}
\newcommand{\efact}{\end{FACT}}
\newcommand{\ehypo}{\end{HYPO}}
\newcommand{\elemm}{\end{LEMM}}
\newcommand{\eprob}{\end{PROB}}
\newcommand{\eprop}{\end{PROP}}
\newcommand{\erema}{\end{REMA}}

\def\0{\mbox{\bf 0}}
\def\1{\mbox{\bf 1}}
\def\2{\mbox{\bf 2}}
\def\3{\mbox{\bf 3}}
\def\4{\mbox{\bf 4}}
\def\5{\mbox{\bf 5}}
\def\6{\mbox{\bf 6}}
\def\7{\mbox{\bf 7}}
\def\8{\mbox{\bf 8}}
\def\9{\mbox{\bf 9}}

\def\b{\mbox{\boldmath $b$}}

\def\u{\mbox{\boldmath $u$}}

\def\x{\mbox{\boldmath $x$}}
\def\y{\mbox{\boldmath $y$}}

\def\A{\mbox{\boldmath $A$}}
\def\B{\mbox{\boldmath $B$}}

\def\O{\mbox{\boldmath $O$}}

\def\U{\mbox{\boldmath $U$}}

\def\W{\mbox{\boldmath $W$}}
\def\X{\mbox{\boldmath $X$}}
\def\Y{\mbox{\boldmath $Y$}}

\def\DC{\mbox{$\cal D$}}
\def\EC{\mbox{$\cal E$}}
\def\FC{\mbox{$\cal F$}}

\def\inprod#1#2{\langle#1, \, #2\rangle}

\def\bdelta{\mbox{\boldmath $\delta$}}

\def\s0{\mbox{\scriptsize \boldmath $0$}}

\def\Real{\mathbb{R}}

\def\SymMat{\mathbb{S}}

\def\rank{{\rm rank}}

\begin{document}

\title{Local-to-Global Exactness of SDP Relaxations\\ for Sparse QCQPs}

\author{
\normalsize
Masakazu Kojima\thanks{Department of Industrial and Systems Engineering,
	Chuo University, Tokyo 192-0393, Japan ({\tt kojima2414@gmail.com}).
 	}, \and \normalsize
Sunyoung Kim\thanks{Department of Mathematics, Ewha W. University, 52 Ewhayeodae-gil, Sudaemoon-gu, Seoul 03760, Korea 
			({\tt skim@ewha.ac.kr}). 
			 }, \and \normalsize
Naohiko Arima\thanks{Independent researcher 
	({\tt nao$\_$arima@me.com}).}
}

\date{\today}

\maketitle 


\begin{abstract}
\noindent
We study exact semidefinite programming (SDP) relaxation for a given
sparse quadratically constrained quadratic program (QCQP). 
The SDP relaxation is exact if, whenever it has an optimal solution, it admits a rank-at-most-one 
optimal solution that corresponds to an optimal solution of the QCQP.
 Using the maximal cliques of a chordal extension of the aggregate sparsity pattern graph
of the data matrices, we formulate the SDP relaxation in terms of
clique-wise matrix variables and develop a local-to-global framework for certifying exactness.
For each clique-wise matrix variable, we introduce a local sub-SDP with
two parameters: a local right-hand-side vector and a consistency matrix
specifying the values of entries shared by overlapping clique-wise matrix 
variables.  In the main theorem, these parameters are determined by an
optimal solution of the global clique-wise SDP.  The theorem shows that
if the resulting local sub-SDPs are exact, 
then the original SDP relaxation is exact.  Under the additional
assumption that any two distinct cliques intersect in at most one node, 
we present three classes of local QCQPs that can be incorporated
into this framework: convex local QCQPs, local QCQPs characterized by sign-pattern 
conditions, and separable local QCQPs 
with a limited number of constraints. 
Examples illustrate how these different local QCQP  
classes can be combined in sparse QCQPs.
\end{abstract}

\noindent
{\bf Key words.}
Quadratically constrained quadratic program, 
semidefinite programming relaxation, 
exact SDP relaxation, 
sparse optimization, 
chordal graph, 
clique-wise formulation, 
local-to-global exactness, 
block-clique sparsity, 
rank-at-most-one optimal solution.

\vspace{0.5cm}

\noindent
{\bf AMS Classification.} 
90C20,  	
90C22,  	
90C25,     
90C26.  	

\section{Introduction}

\label{section:Introduction}

Quadratically constrained quadratic programs (QCQPs) form a broad class
of nonconvex optimization problems.  They include many fundamental models
in operations research, control, signal processing, and combinatorial
optimization. A standard approach to such problems is to lift the quadratic
terms to a symmetric matrix variable and to relax the resulting rank-at-most-one constraint,
thereby obtaining a semidefinite programming (SDP) relaxation.  SDP
relaxations of QCQPs have been extensively studied; see, for example,
\cite{BAO2011,BENTAL2001,FUJIE1997,LUO2010,SHOR1987}
and the references therein.  
A central question is when such an SDP relaxation is exact, in the sense 
that it admits a rank-at-most-one optimal solution, which corresponds to 
an optimal solution of the original QCQP. (We refer to an optimal solution 
of rank at most one as a rank-at-most-one optimal solution.)

The exactness of SDP relaxations has been studied from several different
viewpoints.  One line of work derives exactness conditions directly in
terms of the data matrices, including convexity
\cite{BENTAL2001} and sign-pattern 
conditions
\cite{AZUMA2022,AZUMA2023,KIM2003,SOJOUDI2014}.  
Another line is based on the number of quadratic constraints and theoretical 
rank bounds for SDP solutions \cite{HUANG2010,LUO2010,PATAKI1998}.  
A further line studies
exactness through geometric properties of the feasible region or of
associated cones, such as non-intersecting quadratic constraints
(NIQC) conditions \cite{ARIMA2026,JOYCE2024,YANG2018} and 
the rank-one-generated (ROG) property
\cite{ARGUE2023,ARIMA2023,KIM2020}.  The 
NIQC and ROG viewpoints are closely related, as discussed in \cite{ARIMA2026},
but they are qualitatively different in nature from coefficient-wise or
constraint-count conditions: they depend on how the quadratic
inequalities jointly shape the feasible region.

This paper develops  a local-to-global exactness framework 
 for sparse QCQPs 
that can be verified locally and then assembled into a global exactness certificate.
Many existing exactness results are formulated for a single
QCQP with a specific global structure as mentioned above.  In sparse problems, however,
different parts of the problem may have different structures: one part
may correspond to a convex subproblem, another to a sign-pattern class,
and another to a separable problem with a limited  number of constraints.
Convexity and sign-pattern conditions are particularly well suited to
such a local treatment, because their exactness guarantees are stable
under changes in local right-hand-side values and in the values shared by
overlapping local subproblems.  Rank-bound arguments based on the 
limited  number
of constraints can also be incorporated, although their
applicability may depend on the parameters induced by the global problem.
Since the ROG and NIQC viewpoints are more global and 
geometric in nature, they fall outside the scope of the present work.

Sparsity has long played an important role in semidefinite programming.
For SDPs whose data matrices have a sparse aggregate pattern, chordal
extensions and positive semidefinite matrix completion 
 enable the replacement of a single large positive semidefinite constraint with
 smaller positive semidefinite constraints on maximal cliques
\cite{GRON1984,FUKUDA2003,NAKATA2003}. The resulting clique-wise SDP is
equivalent to the original SDP relaxation, but  it is formulated in terms of 
local matrix variables linked by  consistency constraints on  overlaps.
We adopt this clique-wise formulation 
for the analysis of exactness. 
We emphasize that 
the clique-wise formulation is used here as a theoretical tool for certifying exactness, 
rather than as a computational strategy for solving the sparse SDP relaxation.
Once exactness has been certified,
the SDP relaxation may be solved by any suitable SDP method.  For
computational methods for sparse SDPs based on chordal decomposition, see
\cite{FUKUDA2003,NAKATA2003}.

The main contribution of this paper is the local-to-global exactness
framework for sparse QCQPs.
For a chordal extension of the aggregate 
sparsity pattern graph, we introduce a clique-wise SDP formulation and
local sub-SDPs associated with the maximal cliques.  
Each local sub-SDP
depends on a pair of parameters: a local right-hand-side vector
$\bdelta$ and a consistency matrix $\U$, which specifies the values of 
 entries shared by  local matrix variables. 
 In the main theorem, these
parameters are induced by an optimal solution of the global clique-wise SDP.  
The theorem shows how these induced local certificates can be assembled:
if each local sub-SDP with its induced parameters $(\bdelta,\U)$ admits a
rank-at-most-one optimal solution, 
then the original SDP relaxation admits a rank-at-most-one
optimal solution, and hence is exact.
This `local-to-global exactness' is established for general chordal extensions.

The remaining difficulty is ensuring consistency of local
rank-at-most-one solutions on  overlaps.  For a general chordal
extension, two maximal cliques may share more than one node.  The
resulting consistency constraints 
may then include off-diagonal entries  of local matrix variables, 
and local rank-at-most-one solutions must satisfy nontrivial product 
relations on those overlapping off-diagonal entries. 
These relations are difficult to verify from local exactness
alone.  Therefore, when applying this framework, 
we impose a block-clique assumption: any two maximal
cliques intersect in at most one node.  Under this assumption all
consistency constraints are diagonal, and consistency on overlaps reduces
to matching squared scalar values.  This permits the combination of local
exactness mechanisms without imposing additional off-diagonal rank-one
consistency conditions.

Under the block-clique assumption, the exactness mechanisms discussed above are
recast as local results. 
Some are independent of the induced parameters $(\bdelta,\U)$, while others,
especially those based on rank bounds for separable subproblems, depend on them. 
 We further prove a preservation result showing that certain
dependent inequality constraints can be added without destroying local
exactness.  These results provide the local building blocks for 
the local-to-global exactness certification.

The theoretical significance of the local-to-global exactness framework lies in treating sparsity
as part of the exactness analysis.  The clique-wise formulation 
relates global rank-at-most-one attainment to local rank-at-most-one
attainment.
This provides a route to
proving exactness of SDP relaxations of sparse QCQPs in settings where no single global
exactness criterion applies to the entire problem.  
In this sense,  the contribution is structural: 
it shows how sparsity can serve as 
a mechanism for certifying exact SDP relaxations. 
The examples in Section~\ref{section:example}
illustrate how local exactness certificates can be assigned to different
parts of the sparse structure and then assembled through diagonal consistency
constraints.

We also mention two related lines of work. First, block-clique graph
structures have appeared in doubly nonnegative (DNN) and completely
positive (CPP) reformulations of quadratic optimization problems
\cite{KIM2020a}.  
Although those works concern DNN and
CPP reformulations rather than SDP exactness 
studied here, their use of block-clique  
structures  is closely related to the block-clique assumption mentioned above.

Second, the present framework is complementary to extension results that
preserve exact SDP relaxations under the addition of constraints on a fixed
variable space \cite{KOJIMA2025}.  These results may be viewed as a
vertical extension of a QCQP, whereas 
the local-to-global exactness framework developed
here gives a horizontal extension: exact sub-QCQPs on different variable
subsets are combined through diagonal consistency constraints.  A
preliminary version of this horizontal viewpoint for separable QCQPs
appeared in \cite{KOJIMA2026}; the present paper develops a more general
sparse framework that combines several local exactness mechanisms through
clique-wise decompositions.

The paper is organized as follows.  Section~\ref{section:Preliminaries}
formulates the QCQP and its SDP relaxation, introduces the aggregate 
sparsity pattern graph, and derives the clique-wise formulation based on
a chordal extension of the graph.  Section~\ref{section:framework}
develops the local-to-global exactness framework.  In particular, 
it introduces the local sub-SDPs associated with the clique-wise formulation
developed in Section~\ref{section:Preliminaries},
proves the main theorem, and explains the role of diagonal consistency under the
block-clique assumption.  Section~\ref{section:subQCQPs} establishes
local exactness results for the three classes of sub-QCQPs used in the
framework: convex sub-QCQPs, sub-QCQPs satisfying sign-pattern 
conditions, and separable sub-QCQPs with a limited
number of constraints.  It also presents a preservation result
for dependent inequality constraints.  Section~\ref{section:example}
gives examples illustrating how these local exactness results can be
combined in sparse QCQPs.  Section~\ref{section:conclusion} concludes
the paper.

\section{QCQP, SDP relaxation, and clique-wise reformulation}

\label{section:Preliminaries}

Let $\Real^n$ be the $n$-dimensional Euclidean space of column vectors $\x = (x_1,\ldots,x_n)$, and $\x^T$ the 
transposed row vector of each $\x\in \Real^n$. 
Let $\SymMat^n$ denote the linear space of $n \times n$ symmetric matrices equipped with 
the inner product
$ 
\inprod{\A}{\B} = \sum_{i=1}^n \sum_{j=1}^n [\A]_{ij} [\B]_{ij} \ \text{for } \A, \B \in \SymMat^n,
$ 
and let $\SymMat^n_+$ be the cone of $n \times n$ symmetric positive semidefinite matrices. 
For $\A \in \SymMat^n$, we often write a quadratic form $\x^T\A\x$ in $\x \in \Real^n$ as $\inprod{\A}{\x\x^T}$. 

\subsection{QCQP and its SDP relaxation}

\label{section:QCQPandSDP}

Let $\A_k \in \SymMat^n$ $(k=0,1,\ldots,m)$ and $\b \in \Real^m$. 
We consider the following QCQP:
\begin{eqnarray}
\zeta & = & 
\inf \left\{ \inprod{\A_0}{\x\x^T}  : \x \in \Real^n, \ 
\inprod{\A_k}{\x\x^T} \trianglelefteq_k b_k \ (k=1,\ldots,m)\right\}  
\nonumber \\
& = & 
\inf \left\{ \inprod{\A_0}{\X}  : \X \in \SymMat^n_+, \ \rank(\X) \leq 1, \ 
\inprod{\A_k}{\X} \trianglelefteq_k b_k \ (k=1,\ldots,m)\right\},
\label{eq:QCQP0}
\end{eqnarray}
where $\trianglelefteq_k$ denotes either `$\le$', `$=$', or `$\ge$'. 
The standard SDP relaxation of \eqref{eq:QCQP0}  is given by
\begin{equation}
\eta = 
\inf \left\{ \inprod{\A_0}{\X}  : \X \in \SymMat^n_+, \ 
\inprod{\A_k}{\X} \trianglelefteq_k b_k \ (k=1,\ldots,m)\right\}.  
\label{eq:SDP0}
\end{equation}

The formulation~\eqref{eq:QCQP0} is written in homogeneous quadratic
form.  Linear terms in an inhomogeneous QCQP can be represented by including,
or adding, a normalization constraint
$ 
X_{ii}=1
$ 
for some index $i$.  For a rank-at-most-one matrix $\X=\x\x^T$, this condition
means $x_i=\pm1$.  Since $\x\x^T=(-\x)(-\x)^T$, we may choose the representative
with $x_i=1$.  Then the terms $2[\A_k]_{ij}x_i x_j$ become linear terms in
the remaining variables $x_j$ $(j\ne i,k=0,\ldots,m)$.
The constraint $X_{ii}=1$ fixes only a diagonal entry of the lifted matrix variable $\X$ 
and therefore does not add any edge to the aggregate sparsity pattern graph
defined below.  In contrast, the off-diagonal entries used to represent linear terms are
treated as part of the data  matrices $\A_k$ $(k=0,\ldots,m)$ and are included in the
aggregate sparsity pattern in the same way as the other quadratic
coefficients.

If QCQP~\eqref{eq:QCQP0} is infeasible, we assume that $\zeta=+\infty$. 
Throughout this paper, {\em exactness} of an SDP relaxation is understood in
the rank-attainment sense. More precisely, if the SDP relaxation has an optimal solution,
then it has an optimal solution of rank at most one. 
Such a rank-at-most-one optimal solution is also  optimal for the 
corresponding QCQP, because the QCQP feasible region is precisely the
rank-at-most-one portion of the SDP feasible region.

\subsection{Aggregate sparsity pattern}

To describe the structured sparsity of QCQP~\eqref{eq:QCQP0},
we introduce the {\it aggregate sparsity pattern graph} $G(N,\EC^0)$ of the data  matrices
$\A_k$ $(k=0,1,\ldots,m)$ with $N = \{1,\ldots,n\}$ and 
\[
\EC^0= \left\{
(i,j)\in N\times N : i\neq j,\ [\A_k]_{ij}\neq 0 \ \mbox{for some } k\in\{0,1,\ldots,m\} \right\}, 
\]
where 
$(i,j)$ and $(j,i)$ are identified, since they both represent the same undirected edge 
between nodes $i$ and $j$ 
$(i \ne j)$. 
The sparsity structure can also be encoded by 
{\em the aggregate sparsity pattern matrix}, which is 
an $n \times n$ symmetric symbolic matrix with * at the $(i,j)$th element  for $(i,j) \in \EC^0$ 
and blank elsewhere; * is assigned at $(i,j)$th element if and only if $[\A_k]_{ij} \not= 0$ for 
some $k\in\{0,\ldots,m\}$. 
Figure 1 shows an example of the aggregate sparsity pattern matrix 
and the aggregate sparsity pattern graph $G(N,\EC^0)$. 
This example is used in Section~\ref{section:clique-wiseReformulation} to illustrate
the definitions introduced below  for the clique-wise reformulation of 
QCQP~\eqref{eq:QCQP0} and its SDP relaxation, 
and is also revisited in the subsequent discussion.

\begin{figure}[t!]  
$
\begin{pmatrix} 
1 & *   &   & *  &       &     &  & \\
*  & 2  & *  &   &       &     & *  &  \\
   &  *  & 3 & *  &     &     &   &  \\
*  &     & *  & 4 &  *   &  *  &    &* \\
   &     &    &  *  & 5   &  *  &   &  \\
   &     &    &  *  &  *   & 6  &     &     \\
*  &     &    &     &     &     & 7   &     \\
   &     &    &  *  &       &    &      &  8   \\
\end{pmatrix}
$

\vspace{-45mm}
\mbox{ \ } \hspace{60mm}
\includegraphics[height=35mm]{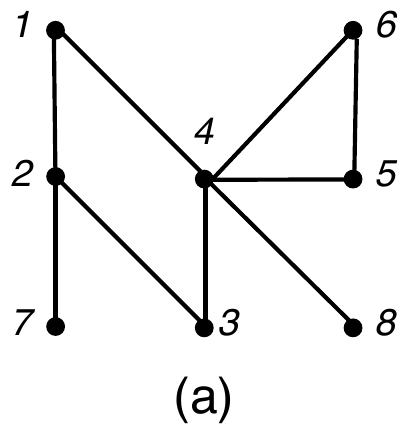}
\mbox{ \ } \hspace{3mm}
\includegraphics[height=35mm]{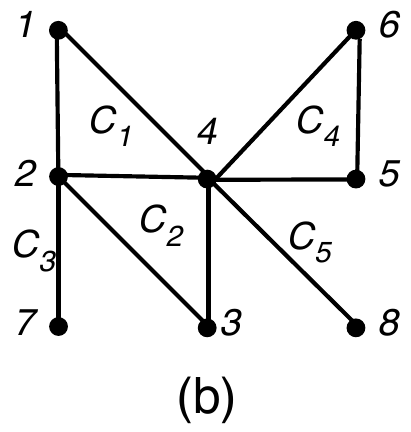}

\vspace{5mm}

\mbox{ \ } \hspace{60mm}
\includegraphics[height=35mm]{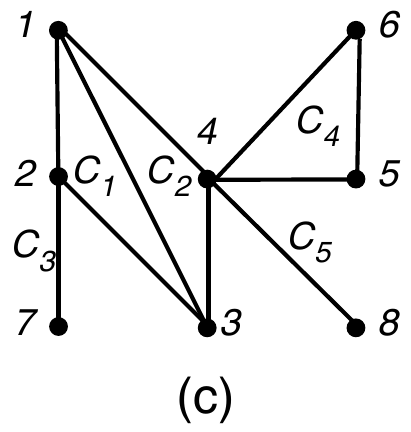}
\mbox{ \ } \hspace{3mm}
\includegraphics[height=35mm]{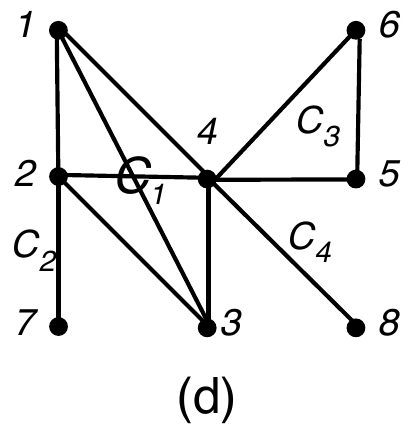}
 
\caption{
An example of the aggregate sparsity pattern matrix (left), 
where * denotes nonzero elements. (a) : the associated aggregate sparsity 
pattern graph $G(N,\EC^0)$ with node set $N = \{1,\ldots,8\}$ and edge set 
$ 
 \EC^0 = \{(1,2),(1,4),(2,3),(2,7),(3,4),(4,5),$ $(4,6),(4,8),(5,6)\}, 
$ 
which is not chordal since the cycle formed by the $4$ edges $(1,2),(2,3),(3,4),(4,1)$ 
is a  chordless 4-cycle. 
(b), (c) and (d): chordal extensions of $G(N,\EC^0)$. In (b), the maximal cliques are 
$C_1=\{1,2,4\},C_2=\{2,3,4\},C_3=\{2,7\},C_4=\{4,5,6\}$, and $C_5=\{4,8\}$. 
In (c),  the maximal cliques are 
$C_1=\{1,2,3\},C_2=\{1,3,4\},C_3=\{2,7\},C_4=\{4,5,6\}$, and $C_5=\{4,8\}$. 
In (d), the maximal cliques are 
$C_1=\{1,2,3,4\},C_2=\{2,7\},C_3=\{4,5,6\}$, and $C_4=\{4,8\}$. 
\label{fig:Example2a}
}
\end{figure} 

%
For every nonempty subset $C$ of $N$, let $\SymMat^C$ denote the linear space of
$|C|\times |C|$ symmetric matrices indexed by $C\times C$, $\SymMat^C_+$ the cone of
positive semidefinite matrices in $\SymMat^C$, and $\Real^C$ the linear space of
$|C|$-dimensional column vectors indexed by $C$.
For every $\A \in \SymMat^n$ and every nonempty subset $C\subseteq N$, we denote by
$\A^C \in \SymMat^C$ the principal submatrix of $\A$ indexed by $C\times C$.
Similarly, for every $\x \in \Real^n$, we denote by $\x^C \in \Real^C$ the subvector
indexed by $C$.

\subsection{A clique-wise reformulation of QCQP~\eqref{eq:QCQP0} and its SDP relaxation}

\label{section:clique-wiseReformulation}

The clique-wise reformulation of QCQP~\eqref{eq:QCQP0} and its SDP
relaxation developed below relies on the following well-known results on
positive semidefinite matrix completions for chordal graphs
\cite{BLAIR1993,GRON1984}.  An undirected graph is said to be
{\em chordal} if every cycle of length at least four has a chord.
For an undirected graph $G(N,\EC)$, define
$\EC \cup \{(i,i): i \in N\}$ by $\overline{\EC}$.
We denote by $\SymMat^n(\EC)$ the class of $n \times n$ partial symmetric matrices $\X$ 
such that
\begin{eqnarray*}
& & [\X]_{ij}=[\X]_{ji}\in \Real \quad \mbox{if } (i,j)\in \overline{\EC}, \
\mbox{and $[\X]_{ij}=[\X]_{ji}$ is left unspecified otherwise}.
\end{eqnarray*}
In particular, if $C$ is a clique of $G(N,\EC)$ and
$\X \in \SymMat^n(\EC)$, then the principal submatrix
$\X^{C} \in \SymMat^C$ is well-defined since
$C \times C \subseteq \overline{\EC}$.

\lemm \label{lemma:chordalGraph}
Let $C_p$ $(p=1,\ldots,\hat{p})$ be the maximal cliques of a chordal
graph $G(N,\EC)$.  Assume that $\X \in \SymMat^n(\EC)$.  Then the
following assertions hold.
\vspace{-2mm}
\begin{description}
\item{(i)}
There exists an $\overline{\X} \in \SymMat^n_+$ such that
$[\X]_{ij} = [\overline{\X}]_{ij}$ for every
$(i,j) \in \overline{\EC}$ if and only if
$\X^{C_p} \in \SymMat^{C_p}_+$ $(p=1,\ldots,\hat{p})$.
We call such an $\overline{\X} \in \SymMat^n_+$ a positive semidefinite
matrix completion of $\X \in \SymMat^n(\EC)$.
\vspace{-2mm}
\item{(ii)}
There exists an $\overline{\X}\in\SymMat^n_+$ such that
$[\X]_{ij}=[\overline{\X}]_{ij}$ for every
$(i,j)\in\overline{\EC}$ and $\rank(\overline{\X})\leq 1$ if and only if
$\X^{C_p}\in\SymMat^{C_p}_+$ and
$\rank(\X^{C_p})\leq 1$ $(p=1,\ldots,\hat{p})$.
We call such an $\overline{\X}$ a rank-at-most-one positive semidefinite
matrix completion of $\X\in\SymMat^n(\EC)$.
\vspace{-2mm}
\end{description}
\elemm
\proof{
Assertion (i) is the positive semidefinite matrix completion theorem for
chordal graphs; see \cite{GRON1984}, and also
\cite{FUKUDA2003,NAKATA2003} for its use in semidefinite programming.
Assertion (ii) follows as a special case of the minimum-rank positive
semidefinite matrix completion theorem for chordal graphs; see
\cite{DANCIS1992,JIANG2023}.
\qed
}

In the remainder of the paper, we let $G(N,\EC)$ be 
a chordal extension of the aggregate
sparsity pattern graph $G(N,\EC^0)$ associated with QCQP~\eqref{eq:QCQP0}.
We denote the maximal cliques of $G(N,\EC)$ by
$C_p$ $(p=1,\ldots,\hat{p})$.  
Let 
$\overline{\EC} = \EC \cup \{(i,i) : i \in N \}$. 
We note that  the aggregate sparsity pattern graph $G(N,\EC^0)$ is uniquely determined 
by the data matrices of 
QCQP~\eqref{eq:QCQP0}, but that there are multiple chordal extensions in general. 
Figure 1 illustrates an example of the aggregate sparsity pattern matrix, 
the aggregate sparsity pattern graph 
$G(N,\EC^0)$, and three different chordal extensions of $G(N,\EC^0)$. 

Since $\EC^0 \subseteq \EC \subseteq \bigcup_{p=1}^{\hat{p}} ( C_p\times C_p)$, we see that 
$[\A_k]_{ij} = 0$ if $(i,j) \not\in \bigcup_{p=1}^{\hat{p}} (C_p\times C_p)$; hence
\begin{eqnarray*}
\inprod{\A_k}{\X} = \sum_{(i,j) \in \bigcup_{p=1}^{\hat{p}}( C_p \times C_p)} [\A_k]_{ij} [\X]_{ij} \ 
\mbox{for every } \X \in \SymMat^n 
\end{eqnarray*}
$(k=0,\ldots,m)$. 
Let $k\in \{0,\ldots,m\}$ be fixed. 
We consider decompositions of the coefficient matrix  $\A_k$, a collection of 
$\A^p_k \in \SymMat^{C_p}$ $(p=1,\ldots,\hat{p})$ such that 
\begin{equation} 
[\A_k]_{ij} = \sum_{p \in P(i,j)} [\A^p_k]_{ij} \ \mbox{for every } (i,j) \ \mbox{such that } P(i,j) \ne \emptyset,  
\label{eq:cliqueDecomposition0}
\end{equation} 
where $P(i,j) = \{ p : (i,j) \in C_p \times C_p\}$. 
Namely, the $(i,j)$th element $[\A_k]_{ij}$ of $\A_k \in \SymMat^n$ is distributed among 
$(C_p \times C_p)$s containing $(i,j)$. 
Then,  
\begin{equation}
\inprod{\A_k}{\X} = \sum_{p=1}^{\hat{p}} \inprod{\A_k^p}{\X^{C_p}}  \ \mbox{for every } \X \in 
\SymMat^n, \label{eq:cliqueDecomposition}
\end{equation} 
or equivalently, 
\[
\sum_{(i,j) \in \bigcup_{p=1}^{\hat{p}}(C_p \times C_p)} [\A_k]_{ij}[\X]_{ij} = \sum_{p=1}^{\hat{p}} \sum_{(i,j) \in C_p\times C_p} [\A_k^p]_{ij}[\X]_{ij}  \ \mbox{for every } \X \in 
\SymMat^n
\]
holds. 
Comparing the terms of symmetric entries $[\X]_{ij}$ with $(i,j) \in (C_p \times C_p)$  
on both sides of the identity, 
we see that the left-hand side contains the term 
$[\A_k]_{ij}[\X]_{ij}$, while the right-hand side contains
$
\sum_{p\in P(i,j)} [\A_k^p]_{ij}[\X]_{ij}.
$
The coefficients of $[\X]_{ij}$ are equal by the  decomposition defined  in 
\eqref{eq:cliqueDecomposition0}. 
Therefore, the identity holds for every  value of $[\X]_{ij}\in\mathbb{R}$.
This implies \eqref{eq:cliqueDecomposition}. 
In the case (b) of Figure 1, we see that 
$ 
P(2,2) = \{1,2,3\}, P(2,4) = \{1,2\}, P(4,4) = \{1,2,4,5\}. 
$ 
Thus, we can take, for example  
\begin{eqnarray*}
& & [\A_k^p]_{44} = \frac{[\A_k]_{44}}{4} \ \mbox{for every } p \in P(4,4) = \{1,2,4,5\} \mbox{ or} \\ 
& &  [\A_k^1]_{44} =  [\A_k^2]_{44} =  [\A_k^4]_{44} = 0, \ [\A_k^5]_{44} = [\A_k]_{44}. 
\end{eqnarray*} 

By substituting identity~\eqref{eq:cliqueDecomposition} into QCQP~\eqref{eq:QCQP0}, 
we obtain an equivalent reformulation of QCQP~\eqref{eq:QCQP0} and 
its SDP relaxation~\eqref{eq:SDP0}.
\begin{eqnarray}
\zeta_{\rm c} 
& = & \inf \left\{ \sum_{p=1}^{\hat{p}}\inprod{\A^p_0}{\X^{C_p}}  : 
\begin{array}{l}
\X \in \SymMat^n(\EC),\X^{C_p} \in \SymMat^{C_p}_+,\rank(\X^{C_p}) \leq 1 \\ 
(p=1,\ldots,\hat{p}),\\ 
\displaystyle \sum_{p=1}^{\hat{p}}\inprod{\A^p_k}{\X^{C_p}} \trianglelefteq_k b_k \ (k=1,\ldots,m)
\end{array}
\right\}.  
\label{eq:QCQP1}
\end{eqnarray}
\begin{eqnarray}
\eta_{\rm c}  = 
\inf \left\{ \sum_{p=1}^{\hat{p}}\inprod{\A^p_0}{\X^{C_p}}  : 
\begin{array}{l}
\X \in \SymMat^n(\EC), \X^{C_p} \in \SymMat^{C_p}_+ \ (p=1,\ldots,\hat{p}), \\ 
\displaystyle  \sum_{p=1}^{\hat{p}}\inprod{\A^p_k}{\X^{C_p}} \trianglelefteq_k b_k \ (k=1,\ldots,m)
\end{array}
\right\}.  
\label{eq:SDP1}
\end{eqnarray}
We say that a feasible solution $\X\in\SymMat^n(\EC)$ of
SDP~\eqref{eq:SDP1} is rank-at-most-one if
\[
\X^{C_p}\in\SymMat^{C_p}_+,\
\rank(\X^{C_p})\leq 1
\ (p=1,\ldots,\hat p). 
\]
Thus a rank-at-most-one feasible solution of SDP~\eqref{eq:SDP1}
is precisely a feasible solution of QCQP~\eqref{eq:QCQP1}; the same
terminology is used for optimal solutions.
By Lemma~\ref{lemma:chordalGraph}(ii), an
$\X\in\SymMat^n(\EC)$ is a feasible solution of
QCQP~\eqref{eq:QCQP1} if and only if it admits a rank-at-most-one
completion $\overline{\X}\in\SymMat^n_+$ that is feasible for the
original QCQP~\eqref{eq:QCQP0}.  The objective values are the same
under this correspondence.  Hence QCQP~\eqref{eq:QCQP0} and
QCQP~\eqref{eq:QCQP1} are equivalent.  Similarly, by
Lemma~\ref{lemma:chordalGraph}(i), the original SDP~\eqref{eq:SDP0}
and SDP~\eqref{eq:SDP1} are equivalent.  Therefore,
$ 
\eta_{\rm c}=\eta\leq \zeta_{\rm c}=\zeta .
$ 
Consequently, the SDP relaxation~\eqref{eq:SDP0} of
QCQP~\eqref{eq:QCQP0} is exact if and only if the SDP
relaxation~\eqref{eq:SDP1} of QCQP~\eqref{eq:QCQP1} is exact, {\it i.e.}, 
the SDP relaxation~\eqref{eq:SDP1} has a rank-at-most-one optimal solution when 
it has an optimal solution. 

Although all constraints in \eqref{eq:QCQP1} and \eqref{eq:SDP1} are expressed
in terms of the clique submatrix variables $\X^{C_p} \in \SymMat^{C_p}$ 
$(p=1,\ldots,\hat{p})$, these variables are not independent in general. 
Indeed, if two sets $C_p$ and $C_q$ overlap, then the corresponding entries of
$\X^{C_p}$ and $\X^{C_q}$ must coincide on 
$(C_p\times C_p)\cap(C_q\times C_q)$.
To describe these consistency requirements explicitly, we introduce local matrix
variables $\Y^p \in \SymMat^{C_p}$ $(p=1,\ldots,\hat{p})$. 
Let
\[
\DC_p=\bigcup_{q\neq p}
\bigl((C_p\cap C_q)\times(C_p\cap C_q)\bigr) \subseteq C_p \times C_p,  
\ (p=1,\ldots,\hat p), 
\]
and define
\[
\overline{\DC}
=
\bigcup_{p=1}^{\hat p}\DC_p
=
\bigcup_{1\leq p<q\leq \hat p}
\bigl((C_p\cap C_q)\times(C_p\cap C_q)\bigr).
\]
Thus, if $\DC_p \ne \emptyset$, then $\DC_p$ is the set of entries of the local matrix variable
$\Y^p$ that are shared
with at least one other local matrix variable,
and $\overline{\DC}$ is the set of all such shared entries. 
Each $\DC_p$ is called {\em a consistency set}. 
The consistency requirement can then be written as
\[ 
[\Y^p]_{ij}=[\Y^q]_{ij}
\
\mbox{whenever }
(i,j)\in (C_p\cap C_q)\times(C_p\cap C_q) \ \mbox{and } p\neq q. 
\] 
Equivalently, we introduce an auxiliary partial symmetric matrix variable
$\U\in\SymMat^n(\overline{\DC})$, called {\em the consistency matrix}, and write
\[
[\Y^p]_{ij}=[\U]_{ij}
\quad ((i,j)\in\DC_p,\ p=1,\ldots,\hat p).
\]
In this formulation, the clique submatrix variables $\X^{C_p}$ in 
\eqref{eq:QCQP1} and \eqref{eq:SDP1} are replaced by local matrix variables 
$\Y^p$. These local matrix variables are linked through the consistency matrix $\U$ 
on their overlapping entries in the consistency sets $\DC_p$ $(p=1,\ldots,\hat p)$.

 We see in case (b) of Figure 1 that 
$\DC_1 =  \DC_2 = \{(2,2),(2,4),(4,4)\}, \DC_3 = \{(2,2)\}, \DC_4 = \DC_5 = \{(4,4)\}$ 
and in case (d) of Figure 1  that $\DC_1 =  \{(2,2),(4,4)\}, \DC_2 = \{(2,2)\}, 
\DC_3 = \DC_4 = \{(4,4)\}$. 
In the former case (b), the consistency 
constraints are 
\begin{eqnarray*}
& & [\Y^1]_{ij} = [\Y^2]_{ij} \ ((i,j) \in \{(2,2),(2,4),(4,4)\}), \  [\Y^1]_{22} = [\Y^3]_{22}, \\
& & [\Y^1]_{44} =[\Y^2]_{44} =[\Y^4]_{44} =[\Y^5]_{44}, 
\end{eqnarray*}
while in the latter case (d), they are 
\begin{eqnarray*}
& & [\Y^1]_{22} = [\Y^2]_{22}, \   [\Y^1]_{44} =[\Y^3]_{44} =[\Y^4]_{44}. 
\end{eqnarray*}
In case (b), the two cliques $C_1$ and $C_2$ are coupled not only through
the diagonal entries $(2,2)$ and $(4,4)$, but also through the
off-diagonal entry $(2,4)$.  In contrast, in case (d), all consistency
constraints are diagonal.

Consequently, SDP~\eqref{eq:SDP1} can be rewritten in the local matrix variables
$\Y^p$ $(p=1,\ldots,\hat p)$ and the consistency matrix $\U$ as follows:
\begin{eqnarray}
\eta_c
=
\inf\left\{
\displaystyle \sum_{p=1}^{\hat{p}} \inprod{\A^p_0}{\Y^p} :
\begin{array}{l}
\U \in \SymMat^n(\overline{\DC}), \
\Y^p \in \SymMat^{C_p}_+ \ (p=1,\ldots,\hat{p}), \\
\displaystyle \sum_{p=1}^{\hat{p}} \inprod{\A^p_k}{\Y^p}
\trianglelefteq_k b_k \ (k=1,\ldots,m), \\
{[\Y^p]_{ij}}=[\U]_{ij} \quad ((i,j)\in \DC_p,\ p=1,\ldots,\hat{p})
\end{array}
\right\}.
\label{eq:SDP4}
\end{eqnarray}
We call~\eqref{eq:SDP4} {\em the clique-wise formulation} of SDP~\eqref{eq:SDP0}, which was originally proposed in \cite{NAKATA2003}.  
The problem obtained from SDP~\eqref{eq:SDP4} by adding a rank-at-most-one condition 
\begin{equation}
\rank(\Y^p)\leq 1 \quad (p=1,\ldots,\hat p)
\label{eq:rankOne}
\end{equation}
is called the clique-wise formulation of QCQP~\eqref{eq:QCQP0}.
Accordingly, a feasible solution $(\Y^1,\ldots,\Y^{\hat p},$ $ \U)$ of
SDP~\eqref{eq:SDP4} is called {\em rank-at-most-one} if it
satisfies~\eqref{eq:rankOne}.
The following two lemmas summarize the relation between the formulations
using the partial matrix $\X$ and those using the local matrix variables 
$\Y^p$ together with the consistency matrix $\U$.

\lemm \label{lemma:mainSDP}
SDPs~\eqref{eq:SDP1} and~\eqref{eq:SDP4} are equivalent.  More precisely,
the following two assertions hold.
\vspace{-2mm}
\begin{description}
\item{(i)}
Suppose that $\X\in \SymMat^n(\EC)$ is a feasible solution of
SDP~\eqref{eq:SDP1} with objective value
$\sum_{p=1}^{\hat p}\inprod{\A^p_0}{\X^{C_p}}$.  Let
\[
\Y^p=\X^{C_p} \ (p=1,\ldots,\hat p), \
\U\in\SymMat^n(\overline{\DC}), \
[\U]_{ij}=[\X]_{ij} \ ((i,j)\in\overline{\DC}).
\]
Then $(\Y^1,\ldots,\Y^{\hat p},\U)$ is a feasible solution of
SDP~\eqref{eq:SDP4} with the same objective value.\vspace{-2mm}
\item{(ii)}
Suppose that $(\Y^1,\ldots,\Y^{\hat p},\U)$ is a feasible solution of
SDP~\eqref{eq:SDP4} with objective value
$\sum_{p=1}^{\hat p}\inprod{\A^p_0}{\Y^p}$.  Then there exists 
a feasible solution 
$\X\in\SymMat^n(\EC)$ of SDP~\eqref{eq:SDP1} satisfying 
$ 
\X^{C_p}=\Y^p \ (p=1,\ldots,\hat p)
$ 
with the same objective  value.
\end{description}
\elemm

\lemm \label{lemma:mainQCQP}
QCQP~\eqref{eq:QCQP1} is equivalent to SDP~\eqref{eq:SDP4} with the
additional rank-at-most-one condition~\eqref{eq:rankOne}.
More precisely, the following two assertions hold.
\vspace{-2mm}
\begin{description}
\item{(i)}
Suppose that $\X\in \SymMat^n(\EC)$ is a feasible solution of
QCQP~\eqref{eq:QCQP1} with  objective value  
$\sum_{p=1}^{\hat p}\inprod{\A^p_0}{\X^{C_p}}$.  Let
\[
\Y^p=\X^{C_p} \ (p=1,\ldots,\hat p), \
\U\in\SymMat^n(\overline{\DC}), \
[\U]_{ij}=[\X]_{ij} \ ((i,j)\in\overline{\DC}).
\]
Then $(\Y^1,\ldots,\Y^{\hat p},\U)$ is a rank-at-most-one feasible solution of the
SDP~\eqref{eq:SDP4} with the same
objective value.
\vspace{-2mm}
\item{(ii)}
Suppose that $(\Y^1,\ldots,\Y^{\hat p},\U)$ is a rank-at-most-one 
feasible solution of SDP~\eqref{eq:SDP4} with objective value 
$\sum_{p=1}^{\hat p}\inprod{\A^p_0}{\Y^p}$.  Then there exists a 
feasible solution $\X\in\SymMat^n(\EC)$ of QCQP~\eqref{eq:QCQP1} satisfying 
$ 
\X^{C_p}=\Y^p \ (p=1,\ldots,\hat p)
$ 
with the same  objective value.
\end{description}
\elemm

By Lemmas~~\ref{lemma:chordalGraph}, \ref{lemma:mainSDP} and~\ref{lemma:mainQCQP}, we have shown that SDP~\eqref{eq:SDP0} and
QCQP~\eqref{eq:QCQP0} are equivalent to their clique-wise formulations, 
SDP~\eqref{eq:SDP4}  and SDP~\eqref{eq:SDP4} with the additional
rank-at-most-one condition~\eqref{eq:rankOne}, respectively. 
Consequently, SDP~\eqref{eq:SDP0} has a rank-at-most-one optimal solution $\X$, 
which is an optimal solution of QCQP~\eqref{eq:QCQP0}, if and only if 
SDP~\eqref{eq:SDP4} has a rank-at-most-one optimal solution 
$(\Y^1,\ldots,\Y^{\hat{p}},\U)$ such that $\X^{C_p} = \Y^p$ $(p=1,\ldots,\hat{p})$. 
This equivalence allows the exactness analysis of the 
original SDP relaxation~\eqref{eq:SDP0} to be carried out through the clique-wise formulation
\eqref{eq:SDP4}.

It should be noted that the choice of decompositions of the coefficient
matrices $\A_k \in \SymMat^n$ into collections
$\A^p_k \in \SymMat^{C_p}$ $(p=1,\ldots,\hat{p})$ satisfying
\eqref{eq:cliqueDecomposition0} does not affect the equivalence of
SDP~\eqref{eq:SDP0} with 
SDP~\eqref{eq:SDP4}.    
The
choice can, however, affect the local sub-SDP~\eqref{eq:subSDP} in  the 
single matrix variable $\Y^p$, 
introduced in the next section, 
where the matrices $\A^p_k$ $(k=0,\ldots,m)$ appear explicitly
in the objective and constraint functions $(p=1,\ldots,\hat{p})$.

\section{Local-to-global exactness via clique subproblems}

\label{section:framework}

Throughout  the remainder of the paper,  SDP~\eqref{eq:SDP4} is referred to as the
global SDP.  For each clique $C_p$ $(p=1,\ldots,\hat p)$, we introduce 
a local sub-SDP, denoted by SDP~\eqref{eq:subSDP}, which describes  the
corresponding local component of the global SDP when the local 
right-hand-side and consistency matrix are fixed. 
The purpose of this section is
to relate the exactness of the global SDP~\eqref{eq:SDP4} to the exactness 
of the associated local sub-SDPs.

\subsection{Clique-wise subproblems}

For $p=1,\ldots,\hat{p}$, we define a local sub-SDP 
in the local matrix variable $\Y^p$ associated with the clique $C_p$ 
by fixing a vector $\bdelta^p = (\delta^p_1,\ldots,\delta^p_m) \in \Real^m$ 
and a consistency matrix $\U \in \SymMat^n(\overline{\DC})$ as parameters.
The parameter $\delta_k^p$ 
represents the distribution of each global 
constraint right-hand side constant $b_k$ into local constraint right-hand side constants 
of the sub-SDPs $(k=1,\ldots,m)$, 
while the other parameter, consistency matrix $\U$ specifies 
the values of the entries of the local matrix variable $\Y^p$ shared by different cliques, {\it i.e.,} 
it fixes the entries of $\Y^p$ on $\DC_p$ so as to enforce consistency 
among overlapping cliques.
We now define the corresponding local sub-SDP  as follows:
\begin{equation}
\eta^p_c(\bdelta^p,\U) = \inf \left\{
\langle \A^p_0, \Y^p \rangle :
\begin{array}{l}
\Y^p \in \SymMat^{C_p}_+, \\ 
\langle \A^p_k, \Y^p \rangle \trianglelefteq_k \delta^p_k \ (k=1,\ldots,m),\\[3pt]
{[\Y^p]_{ij}} = [\U]_{ij} \ ((i,j)\in \DC_p)
\end{array}
\right\}. \label{eq:subSDP}
\end{equation}
The sub-SDP~\eqref{eq:subSDP} is  the local problem associated with $C_p$,
obtained from  the global SDP~\eqref{eq:SDP4}  by fixing  parameters $(\bdelta^p,\U)$.
It will be used to characterize the exactness of SDP~\eqref{eq:SDP4}. 
If  the rank-at-most-one condition rank($\Y^p) \leq 1$ is added to sub-SDP~\eqref{eq:subSDP}, 
then we obtain a local sub-QCQP. 

In the example illustrated in Figure~1(d),
there are four local sub-SDPs associated with
$C_1,C_2,C_3,C_4$.
These sub-SDPs are coupled only through the consistency
matrix $\U \in \SymMat^n(\overline{\DC})$
via the consistency constraints 
$[\Y^p]_{ij}=[\U]_{ij}$ for $(i,j)\in\DC_p$ $(p=1,\ldots,4)$, where
$\DC_1 = \{(2,2),(4,4)\}, \DC_2 = \{(2,2)\}, \DC_3 = \{(4,4)\}, \DC_4 = \{(4,4)\}. $ 

The connection between the local sub-SDPs
\eqref{eq:subSDP} $(p=1,\ldots,\hat p)$ and the global
SDP~\eqref{eq:SDP4} can be expressed through a bilevel optimization problem.
In this representation, the upper level determines the local right-hand-side
vectors $\bdelta^p$ $(p=1,\ldots,\hat p)$, and the common consistency matrix
$\U$, while the lower level consists of the local sub-SDPs
\eqref{eq:subSDP}.  The upper-level objective is the sum of the lower-level
optimal values $\eta^p_c(\bdelta^p,\U)$ $(p=1,\ldots,\hat{p})$.
Define the upper-level problem associated with the global SDP by
\begin{eqnarray}
\widetilde{\eta}_c
&=&
\inf\left\{
\sum_{p=1}^{\hat{p}} \eta^p_c(\bdelta^p,\U) :
\begin{array}{l}
\eta^p_c(\bdelta^p,\U) < +\infty
\quad (p=1,\ldots,\hat{p}), \\
\displaystyle \sum_{p=1}^{\hat{p}} \delta^p_k \trianglelefteq_k b_k
\quad (k=1,\ldots,m), \\
\U \in \SymMat^n(\overline{\DC})
\end{array}
\right\}.
\label{eq:bilevelSDP}
\end{eqnarray}
The conditions
$\eta^p_c(\bdelta^p,\U)<+\infty$ $(p=1,\ldots,\hat{p})$ exclude infeasible local subproblems.
However, they do not exclude   local subproblems  with optimal value $-\infty$ at this stage. 
 
\lemm \label{lemma:valueFunction}
The problem~\eqref{eq:bilevelSDP} has the same optimal value as
SDP~\eqref{eq:SDP4}; that is,
$ 
\widetilde{\eta}_c=\eta_c .
$ 
\vspace{-2mm}
\elemm
\proof{
First, let $(\Y^1,\ldots,\Y^{\hat p},\U)$ be any feasible solution of
SDP~\eqref{eq:SDP4}.  Define
\[
\delta_k^p=\inprod{\A_k^p}{\Y^p}
\quad
(k=1,\ldots,m,\ p=1,\ldots,\hat p).
\]
Then $(\bdelta^1,\ldots,\bdelta^{\hat p},\U)$ is feasible for
problem~\eqref{eq:bilevelSDP}.  Since each $\Y^p$ is feasible for the
corresponding sub-SDP~\eqref{eq:subSDP},
\[
\widetilde{\eta}_c
\leq
\sum_{p=1}^{\hat p}\eta_c^p(\bdelta^p,\U)
\leq
\sum_{p=1}^{\hat p}\inprod{\A_0^p}{\Y^p}.
\]
Taking the infimum over all feasible
$(\Y^1,\ldots,\Y^{\hat p},\U)$ gives
$\widetilde{\eta}_c\leq \eta_c$.

Conversely, fix any feasible solution
$(\bdelta^1,\ldots,\bdelta^{\hat p},\U)$ of
problem~\eqref{eq:bilevelSDP}.  For any feasible solutions
$\Y^p$ of the corresponding sub-SDPs~\eqref{eq:subSDP}
$(p=1,\ldots,\hat p)$, the tuple
$(\Y^1,\ldots,\Y^{\hat p},\U)$ is feasible for SDP~\eqref{eq:SDP4}.
Therefore
$ 
\eta_c
\leq
\sum_{p=1}^{\hat p}\inprod{\A_0^p}{\Y^p}.
$ 
Taking the infimum over the feasible sets of the local sub-SDPs yields
$ 
\eta_c
\leq
\sum_{p=1}^{\hat p}\eta_c^p(\bdelta^p,\U).
$ 
Finally, taking the infimum over all feasible
$(\bdelta^1,\ldots,\bdelta^{\hat p},\U)$ in
problem~\eqref{eq:bilevelSDP} gives
$\eta_c\leq\widetilde{\eta}_c$.  Hence
$\widetilde{\eta}_c=\eta_c$.
\qed
}

\subsection{Relation between the global SDP and local sub-SDPs}

The following theorem is the main result of the paper.  It provides a
local-to-global exactness framework for the clique-wise formulation.
Specifically, an optimal solution
$(\widetilde{\Y}^1,\ldots,\widetilde{\Y}^{\hat p},\widetilde{\U})$
of the global clique-wise SDP~\eqref{eq:SDP4} induces the local right-hand-side vectors
$\tilde{\bdelta}^p$ $(p=1,\ldots,\hat p)$ and the consistency matrix
$\widetilde{\U}$.  The relevant local sub-SDPs~\eqref{eq:subSDP} 
are precisely those
defined by these induced parameters.  If these induced local sub-SDPs 
admit rank-at-most-one optimal solutions, then exactness of the original
SDP relaxation follows.
\theo \label{theorem:main1}
Let $(\widetilde{\Y}^1,\ldots,\widetilde{\Y}^{\hat{p}},\widetilde{ \U})$ be an optimal solution of SDP~\eqref{eq:SDP4}. 
Define 
\begin{eqnarray}
 & & \tilde{\bdelta}^p \in \Real^m, \ 
\tilde{\delta}^p_k = \inprod{\A^p_k}{\widetilde{\Y}^p} 
\ (k=1,\ldots,m, \ p=1,\ldots,\hat{p}). \label{eq:bdelta} 
\end{eqnarray}
Then the following assertions hold:
\vspace{-2mm}
\begin{description}
\item{(i)} 
Each $\widetilde{\Y}^p$ is an optimal solution of sub-SDP~\eqref{eq:subSDP} 
with the parameters  $(\bdelta^p,\U) = (\tilde{\bdelta}^p,\widetilde{\U})$; hence 
$\inprod{\A^p_0}{\widetilde{\Y}^p} = \eta^p_c(\tilde{\bdelta}^p,\widetilde{\U})$ $(p=1,\ldots,\hat p)$ and 
$\sum_{p=1}^{\hat{p}} \eta^p_c(\tilde{\bdelta}^p,\widetilde{\U}) = \sum_{p=1}^{\hat{p}} \inprod{\A^p_0}{\widetilde{\Y}^p} = \eta_c$. \vspace{-2mm}
\item{(ii)} 
For $p =1,\ldots,\hat{p}$, let 
$\widehat{\Y}^p \in \SymMat^{C_p}_+$ be an optimal solution of 
sub-SDP~\eqref{eq:subSDP} with the parameters $(\bdelta^p,\U) = (\tilde{\bdelta}^p,\widetilde{\U})$; 
hence $\inprod{\A^p_0}{\widehat{\Y}^p} = \eta^p_c(\tilde{\bdelta}^p,\widetilde{\U})$. 
Then $(\widehat{\Y}^1,\ldots,\widehat{\Y}^{\hat{p}},\widetilde{ \U})$ 
is an optimal solution of SDP~\eqref{eq:SDP4}, and therefore 
$\eta_c = \sum_{p=1}^{\hat{p}}\inprod{\A^p_0}{\widehat{\Y}^p}$. 
\vspace{-2mm}
\item{(iii) } 
Assume that, for $p =1,\ldots,\hat{p}$, sub-SDP~\eqref{eq:subSDP} with the parameters
$(\bdelta^p,\U) = (\tilde{\bdelta}^p,\widetilde{\U})$ is exact. Then SDP~\eqref{eq:SDP4} is exact. 
\vspace{-2mm}
\end{description}
\etheo
\proof{
(i) 
Let $p \in \{1,\ldots,\hat{p}\}$ be fixed. Then $\widetilde{\Y}^p \in \SymMat^{C_p}_+$, 
$\inprod{\A^p_k}{\widetilde{\Y}^p}  \trianglelefteq_k \tilde{\delta}^p_k,$
$(k=1,\ldots,m)$ and
$[\widetilde{\Y}^p]_{ij}=[\widetilde{\U}]_{ij} \ ((i,j) \in \DC_p)$ 
are satisfied. 
Hence $\widetilde{\Y}^p$ is feasible for sub-SDP~\eqref{eq:subSDP}
with the parameters  $(\bdelta^p,\U) = (\tilde{\bdelta}^p,\widetilde{\U})$.
Assume to the contrary that $\widetilde{\Y}^p$ is not optimal.
Then there exists a feasible solution $\Y^p \in \SymMat^{C_p}_+$ of
sub-SDP~\eqref{eq:subSDP} with the parameters
$(\bdelta^p,\U) = (\tilde{\bdelta}^p,\widetilde{\U})$ 
such that
\[
\inprod{\A^p_0}{\Y^p}< \inprod{\A^p_0}{\widetilde{\Y}^p}, \
 \inprod{\A^p_k}{\Y^p}  \trianglelefteq_k \tilde{\delta}^p_k \ (k=1,\ldots,m), \ 
 [\Y^p]_{ij}= [\widetilde{\U}]_{ij} \ ((i,j)  \in \DC_p). 
\]
Replacing only the $p$th block $\widetilde{\Y}^p$ in the feasible solution 
$(\widetilde{\Y}^1,\ldots,\widetilde{\Y}^{\hat{p}},\widetilde{\U})$ of SDP~\eqref{eq:SDP4} 
by $\Y^p$, we obtain a feasible solution of SDP~\eqref{eq:SDP4}, 
which we denote as $(\overline{\Y}^1,\ldots,\overline{\Y}^{\hat{p}},\widetilde{\U})$, 
where $\overline{\Y}^p = \Y^p$ and $\overline{\Y}^q = \widetilde{\Y}^q$ $(q \ne p)$. 
Then 
\begin{eqnarray*}
\sum_{q=1}^{\hat{p}}  \inprod{\A^q_0} {\overline{\Y}^q} = 
\sum_{q\ne p}  \inprod{\A^q_0} {\widetilde{\Y}^q} 
+ \inprod{\A^p_0} {\overline{\Y}^p} 
< 
\sum_{q=1}^{\hat{p}}  \inprod{\A^q_0} {\widetilde{\Y}^q}
\end{eqnarray*}
holds. This contradicts the optimality of $(\widetilde{\Y}^1,\ldots,\widetilde{\Y}^{\hat{p}},\widetilde{\U})$ 
for SDP~\eqref{eq:SDP4}. 
Therefore $\widetilde{\Y}^p$ is an optimal solution of
sub-SDP~\eqref{eq:subSDP} with the parameters  $(\bdelta^p,\U) = (\tilde{\bdelta}^p,\widetilde{\U})$.

(ii) 
By feasibility of each $\widehat{\Y}^p \in \SymMat^{C_p}$ for sub-SDP~\eqref{eq:subSDP},
$
\inprod{\A^p_k}{\widehat{\Y}^p}
\ \trianglelefteq_k\ \delta_k^p
\ (k=1,\ldots,m,\ p=1,\ldots,\hat{p}). 
$
Summing the $k$th inequalities over $p=1,\ldots,\hat{p}$, we obtain
\[
\sum_{p=1}^{\hat{p}}\inprod{\A^p_k}{\widehat{\Y}^p}
\ \trianglelefteq_k\
\sum_{p=1}^{\hat{p}}\delta_k^p
=
\sum_{p=1}^{\hat{p}}\inprod{\A^p_k}{\widetilde{\Y}^p} 
 \trianglelefteq_k\ b_k
\quad (k=1,\ldots,m).
\]
We also see that the other consistency constraints $[\widehat{\Y}^p]_{ij}=[\widetilde{\U}]_{ij} 
\ ((i,j)\in \DC_p,p=1,\ldots,\hat{p})$ are satisfied by the feasibility of each 
$\widehat{\Y}^p$ for sub-SDP~\eqref{eq:subSDP} $(p=1,\ldots,\hat{p})$.
Thus $(\widehat{\Y}^1,\ldots,\widehat{\Y}^{\hat{p}},\widetilde{\U})$ is feasible for
SDP~\eqref{eq:SDP4}.

Next, we show its optimality.
By Assertion (i), for each $p=1,\ldots,\hat{p}$,
$\widetilde{\Y}^p$ is an optimal solution of the same subproblem
sub-SDP~\eqref{eq:subSDP} with the parameters  $(\bdelta^p,\U) = (\tilde{\bdelta}^p,\widetilde{\U})$.
Since $\widehat{\Y}^p$ is also optimal for that subproblem, we have
$
\inprod{\A^p_0}{\widehat{\Y}^p}
=
\inprod{\A^p_0}{\widetilde{\Y}^p}
\ (p=1,\ldots,\hat{p}).
$
Summing these equalities over $p$ yields
$
\sum_{p=1}^{\hat{p}}\inprod{\A^p_0}{\widehat{\Y}^p}
=
\sum_{p=1}^{\hat{p}}\inprod{\A^p_0}{\widetilde{\Y}^p}.
$
The right-hand side is the optimal value of SDP~\eqref{eq:SDP4}, because
$(\widetilde{\Y}^1,\ldots,\widetilde{\Y}^{\hat{p}},\widetilde{\U})$ is optimal for
SDP~\eqref{eq:SDP4}. Therefore,
$(\widehat{\Y}^1,\ldots,\widehat{\Y}^{\hat{p}},\widetilde{\U})$ is also an optimal solution 
of SDP~\eqref{eq:SDP4}. 

(iii) 
By assumption, $(\widetilde{\Y}^1,\ldots,\widetilde{\Y}^{\hat p},\widetilde{\U})$ is an 
optimal solution of  SDP~\eqref{eq:SDP4}.
It suffices to show that  SDP~\eqref{eq:SDP4} has a 
rank-at-most-one optimal solution. 
By Assertion~(i), each $\widetilde{\Y}^p$ is an optimal solution of
sub-SDP~\eqref{eq:subSDP} with the parameters
$(\bdelta^p,\U)=(\tilde{\bdelta}^p,\widetilde{\U})$.
By the exactness assumption on this sub-SDP, 
there exists a rank-at-most-one optimal solution $\widehat{\Y}^p$ of this
sub-SDP. 
 By Assertion~(ii),
$(\widehat{\Y}^1,\ldots,\widehat{\Y}^{\hat p},\widetilde{\U})$ is an
optimal solution of SDP~\eqref{eq:SDP4}.  Since
$\rank(\widehat{\Y}^p)\leq 1$ $(p=1,\ldots,\hat{p})$, this
optimal solution is rank-at-most-one.  Hence, SDP~\eqref{eq:SDP4} is
exact.
\qed
}

\subsection{Diagonal consistency assumption}

\label{section:diagonal}

Theorem~\ref{theorem:main1}(iii) requires exactness of the local
sub-SDPs~\eqref{eq:subSDP} with the parameters
$(\bdelta^p,\U)=(\tilde{\bdelta}^p,\widetilde{\U})$ induced by an
optimal solution of the global SDP~\eqref{eq:SDP4}
$(p=1,\ldots,\hat p)$.  
These parameters are not known a priori.
Therefore, in order to use the theorem together with local exactness results, 
it is necessary to identify classes of local sub-SDPs whose exactness can be guaranteed for
all admissible parameter values, or under conditions that can be checked
without solving the global SDP.

A critical obstruction for identifying such classes of sub-SDPs 
is caused by the unknown consistency matrix $\widetilde{\U}$.
To see this, focus only on the consistency constraints in the 
sub-SDP~\eqref{eq:subSDP} with the parameter 
$(\bdelta^p,\U)=(\tilde{\bdelta}^p,\widetilde{\U})$,
\[
[\Y^p]_{ij}=[\widetilde{\U}]_{ij}
\quad ((i,j)\in\DC_p).
\]
For these constraints to be satisfied by an unknown rank-at-most-one optimal solution 
$\widetilde{\Y}^p$ of the sub-SDP~\eqref{eq:subSDP}, 
the entries of $\widetilde{\U}$ must satisfy the necessary condition
\begin{equation}
[\u]_i[\u]_j = [\widetilde{\U}]_{ij}
\quad ((i,j)\in\DC_p)
\quad \mbox{for some } \u \in \Real^{C_p}.
\label{eq:compRelation}
\end{equation}
Thus, when the consistency set $\DC_p$ contains off-diagonal entries, the 
consistency matrix $\widetilde{\U}$ must satisfy the
nontrivial product relations~\eqref{eq:compRelation}.  
If \eqref{eq:compRelation} is not satisfied, the sub-SDP~\eqref{eq:subSDP} 
cannot have a  rank-at-most-one feasible solution. 
Also, \eqref{eq:compRelation} cannot be verified in advance 
since the consistency matrix $\widetilde{\U}$ is unknown. 

To avoid this obstruction, we impose the following diagonal consistency
assumption:
\begin{equation}
\DC_p \subseteq \{(i,i): i\in N\}
\quad (p=1,\ldots,\hat p).
\label{eq:blockClique}
\end{equation}
Equivalently, any two distinct maximal cliques intersect in at most one
node,
\[
|C_p\cap C_q|\leq 1 \quad (p\ne q).
\]
A chordal graph satisfying this assumption is called a block-clique graph
\cite{JOHNSON1996,KIM2020a}.  
Under the diagonal consistency assumption, \eqref{eq:compRelation}
reduces to the simpler necessary condition 
\[
([\u]_i)^2 = [\widetilde{\U}]_{ii}
\quad ((i,i)\in\DC_p)
\quad \mbox{for some } \u \in \Real^{C_p}.
\]
Since $[\widetilde{\U}]_{ii}=[\widetilde{\Y}^p]_{ii}\geq0$ for every
$(i,i)\in\DC_p$, this condition is always satisfied, and the 
obstruction~\eqref{eq:compRelation} 
on the unknown consistency matrix $\widetilde{\U}$ caused 
by off-diagonal entries in $\DC_p$ has been removed. 

It is important to note that
the block-clique assumption does not imply exactness of the 
sub-SDP~\eqref{eq:subSDP}.  Its role is to remove the 
off-diagonal product obstruction on $\widetilde{\U}$ caused by
\eqref{eq:compRelation} when applying Theorem~\ref{theorem:main1}(iii).
Thus the block-clique assumption makes it possible to formulate local
exactness results in Section~\ref{section:subQCQPs} without having to
verify unknown off-diagonal product relations \eqref{eq:compRelation} 
on the consistency matrix $\widetilde{\U}$.

The diagonal consistency assumption~\eqref{eq:blockClique} leads to the following 
three representative cases of the consistency set $\DC_p$:
\vspace{-2mm}
\begin{description}
\item{(0)} $\DC_p=\emptyset$.
\vspace{-2mm}
\item{(I)} $\DC_p=\{(i,i)\}$ for some $i\in C_p$.
\vspace{-2mm}
\item{(II)} $\{(i,i):i\in C_p\} \supseteq \DC_p$ with multiple diagonal indices.
\vspace{-2mm}
\end{description}
In case (0), the local matrix variable $\Y^p$ is not subject to any consistency constraints. 
The corresponding subproblem is therefore independent of the other
local sub-SDPs at the level of matrix variables, although coupling may still 
occur  through the right-hand-side parameter vector $\bdelta^p$.

In case (I), the consistency constraints fix one diagonal entry of the local matrix variable
$\Y^p$.  This occurs when the clique $C_p$ shares a single node with another clique,
as illustrated in the cases of $\DC_2$, $\DC_3$, and $\DC_4$ in Figure~1(d), where
$[\Y^2]_{22}=[\U]_{22}$, $[\Y^3]_{44}=[\U]_{44}$, and
$[\Y^4]_{44}=[\U]_{44}$ are consistency constraints, respectively.

In case (II), multiple diagonal entries of the local matrix variable $\Y^p$ are fixed through the
consistency matrix $\U$.  This may occur when $C_p$ intersects several
other cliques, each in a single node. 
This is illustrated by  the
case of $C_1$ in Figure~1(d), where the consistency constraints are $[\Y^1]_{22}=[\U]_{22}$ and
$[\Y^1]_{44}=[\U]_{44}$. 

In the next section, we identify three classes of local sub-QCQPs whose
SDP relaxations admit rank-at-most-one optimal solutions under this
diagonal consistency assumption~\eqref{eq:blockClique}.  These classes provide local building
blocks for the local-to-global exactness framework.

\section{Sub-QCQPs with exact SDP relaxations}

\label{section:subQCQPs}

In this section, we present three classes of local sub-QCQPs whose 
corresponding local sub-SDPs are exact. These  are: 
convex sub-QCQPs,
sub-QCQPs characterized by sign-pattern 
conditions, and
separable sub-QCQPs with a limited number of constraints.
Sections~\ref{section:convex}, \ref{section:signPattern}, and~\ref{section:separable} 
discuss these three classes, respectively.
For example, in Figure~1(d), different classes can be assigned to
different clique-wise subproblems.  Further examples will be presented
in Section~\ref{section:example}.

The first two classes have a parameter-independent character: their
local SDP relaxations are exact for every 
right-hand-side vector and consistency matrix. 
By contrast, the 
exactness of the third class may depend on the right-hand-side vector 
induced from an optimal solution of the global SDP.
This distinction is important when the local exactness results are combined
through Theorem~\ref{theorem:main1}(iii).

For notational simplicity, 
we fix a clique
$C_p$ throughout this section and suppress the superscript $p$.  Thus we write
$C_p=C=\{1,\ldots,\ell\}$,
$\A^p_k=\A_k$, $\bdelta^p=\bdelta$, and
$\DC_p=\DC$.  
Under the diagonal consistency assumption~\eqref{eq:blockClique}, we have
$\DC\subseteq\{(i,i):1\le i\le \ell\}$.
Then the local sub-SDP is written as
\begin{equation}
\eta(\bdelta,\U) = \inf \left\{
\langle \A_0,\Y \rangle :
\begin{array}{l}
\Y \in \SymMat^{\ell}_+, \\
\langle \A_k, \Y \rangle \trianglelefteq_k \delta_k \ (k=1,\ldots,m),\\[3pt]
{[\Y]_{ii}} = U_{ii} \ ((i,i)\in \DC)
\end{array}
\right\}.
\label{eq:subSDP1}
\end{equation}
We call the problem obtained by adding
$\rank(\Y)\leq1$ to SDP~\eqref{eq:subSDP1} the local sub-QCQP.  

\subsection{Sub-QCQPs characterized by convexity} 

\label{section:convex}

\theo \label{theorem:convex}
Assume that `$\trianglelefteq_k$' $=$ `$\le$' $(1 \leq k \leq m)$.
Suppose that one of the following conditions holds:
\vspace{-2mm}
\begin{description}
\item{(0)} $\DC = \emptyset$ and $\A_k$ $(k=0,\ldots,m)$ are positive semidefinite.\vspace{-2mm}
\item{(I)} $\DC = \{(i,i)\}$ for some $i \in C$ and $\A_k^{C\backslash\{i\}}$
$(k=0,\ldots,m)$ are positive semidefinite.\vspace{-2mm}
\end{description}
Then, sub-SDP~\eqref{eq:subSDP1} is exact, {\it i.e.},
sub-SDP~\eqref{eq:subSDP1} has a rank-at-most-one optimal solution
whenever it has an optimal solution, for every $\bdelta\in\Real^m$ and $\U \in \SymMat^{\ell}(\DC)$.\vspace{-2mm}
\etheo
\proof{
Let $\bdelta\in\Real^m$ and $\U \in \SymMat^{\ell}(\DC)$ be fixed arbitrarily. 
If sub-SDP~\eqref{eq:subSDP1} has no optimal solution, there is nothing
to prove.  Suppose that it has an optimal solution.  In case (I), 
$U_{ii}\geq0$ $((i,i) \in \DC)$ by positive semidefiniteness of $\Y$. 
The QCQP obtained from sub-SDP~\eqref{eq:subSDP1} with adding 
the rank constraint rank$(\Y) \leq 1$ can be written as 
\begin{equation}
\zeta(\bdelta,\U) = \inf \left\{
\langle \A_0,\y\y^T \rangle :
\begin{array}{l}
\langle \A_k, \y\y^T\rangle \leq \delta_k \quad (k=1,\ldots,m),\\[3pt]
y_i = +\sqrt{U_{ii}} \ \mbox{or } y_i = -\sqrt{U_{ii}}
\quad ((i,i)\in \DC)
\end{array}
\right\}.
\label{eq:subQCQP2}
\end{equation}
In case (0), $\DC=\emptyset$, and the quadratic function
$\inprod{\A_k}{\y\y^T}$ is convex in $\y\in\Real^\ell$ 
$(k=0,\ldots,m)$.

In case (I), for simplicity of notation, assume $i=\ell$.  Since
$\inprod{\A_k}{\y\y^T}=\inprod{\A_k}{(-\y)(-\y)^T}$ for every
$k=0,\ldots,m$, the two choices
$y_\ell=+\sqrt{U_{\ell\ell}}$ and
$y_\ell=-\sqrt{U_{\ell\ell}}$ are equivalent up to the global sign change
$\y\mapsto-\y$.  Hence, we may fix
$y_\ell=+\sqrt{U_{\ell\ell}}$.  After this fixing, the quadratic
function $\inprod{\A_k}{\y\y^T}$ is convex in the remaining variables
$(y_1,\ldots,y_{\ell-1})$ since
$\A_k^{C\backslash\{\ell\}}\succeq0$ $(k=0,\ldots,m)$.

Thus, in both cases, sub-SDP~\eqref{eq:subSDP1} can be viewed as the SDP
relaxation of a convex QCQP.  It is well-known that the SDP relaxation is
exact for such convex QCQPs. See, for example, \cite{SHOR1987} and
\cite[Section 4.2]{BENTAL2001}.  Therefore,
sub-SDP~\eqref{eq:subSDP1} has a rank-at-most-one optimal solution
whenever it has an optimal solution.
\qed
}

We note that, except for the assumptions in (0) and (I), no
additional assumption is required on the right-hand-side vector $\bdelta$ or the
consistency matrix $\U$.

\rema
The argument above does not extend directly to the case $|\DC|\geq2$.
For example, suppose that
$\DC=\{(1,1),(2,2)\}$.  Fixing
$y_1=+\sqrt{U_{11}}$ and $y_2=+\sqrt{U_{22}}$ is not, in general,
equivalent to fixing
$y_1=+\sqrt{U_{11}}$ and $y_2=-\sqrt{U_{22}}$. 
Thus the rank-constrained
problem~\eqref{eq:subQCQP2}  cannot be reduced consistently to a single convex QCQP by fixing
the signs of the variables indexed by $\DC$.  This is why the convex
case in Theorem~\ref{theorem:convex} is restricted to
$\DC=\emptyset$ or $|\DC|=1$.
\erema

\subsection{Sub-QCQPs characterized by sign-pattern conditions} 

\label{section:signPattern}

We next  consider a class of generally nonconvex sub-QCQPs whose SDP relaxations 
are exact under suitable sign-pattern conditions on the data matrices $\A_k$ $(0\leq k \leq m)$. 
Throughout this section, we assume that $\trianglelefteq_k$ $=$ `$\le$' $(1 \leq k \leq m)$. 
As in the convex case, the exactness result for this class holds for every 
right-hand-side vector $\bdelta$ and consistency matrix $\U$.
Let $G(L,\FC)$ denote the aggregate sparsity pattern graph 
with the node set $L = \{1,\ldots,\ell\}$ and 
\begin{eqnarray*}
\FC & = & \left\{ (i,j) \in L \times L : i \not= j, \ [\A_k]_{ij} \ne 0 \ \mbox{for some } 
k \in \{0,1,\ldots,m\} \right\}.  
\end{eqnarray*}
We assume that $\DC \subseteq \{(i,i) : i  \in C\}$ and $0 \leq |\DC| \leq \ell$  (case (II)). 
Note that the equality constraints 
$[\Y]_{ii} = U_{ii}$ $((i,i) \in \DC)$ in sub-SDP~\eqref{eq:subSDP1} and the 
corresponding sub-QCQP can be replaced by the inequality constraints 
$ 
[\Y]_{ii} \leq U_{ii}, \ -[\Y]_{ii} \leq -U_{ii} \ ((i,i) \in \DC). 
$ 
These inequality constraints, however, do not affect the aggregate
sparsity pattern graph $G(L,\FC)$ of their data matrices.

For every $(i,j) \in \FC$, define 
\begin{eqnarray*}
\sigma_{ij} =
\begin{cases}
+1 & \text{if } [\A_k]_{ij} \ge 0 \text{ for all } k \in \{0,1,\ldots,m\},\\
-1 & \text{if } [\A_k]_{ij} \le 0 \text{ for all } k \in \{0,1,\ldots,m\},\\
0  & \text{otherwise}.
\end{cases}
\end{eqnarray*}
Let $\{F_1,\ldots,F_r\}$ denote a cycle basis for $G(L,\FC)$. The following theorem and its corollary 
follow directly from \cite[Theorem 2 and Corollary 1]{SOJOUDI2014},  
respectively. 
\theo \label{theorem:signDefinite} (\cite[Theorem 2]{SOJOUDI2014}) Assume that 
\vspace{-2mm} 
\begin{description}
\item{(i) } $\sigma_{ij} \in \{-1,1\} \ \mbox{for every } (i,j) \in \FC$, \vspace{-1mm} 
\item{(ii) } $\prod_{(i,j)\in F_s} \sigma_{ij} = (-1)^{|F_s|} \mbox{for every } s=1,\ldots,r$. \vspace{-1mm} 
\end{description}
Then sub-SDP~\eqref{eq:subSDP1}  is exact for every $\bdelta\in\Real^m$ and $\U \in \SymMat^{\ell}(\DC)$. 
\etheo
\coro \label{corollary:signDefinite} (\cite[Corollary 1]{SOJOUDI2014})
Assume that one of the following conditions holds: \vspace{-2mm}
\begin{description}
\item{(i) } The graph $G(L,\FC)$ is arbitrary and $\sigma_{ij} = -1$ 
for every $(i,j) \in \FC$ 
(or equivalently all off-diagonal elements of $\A_k$  are nonpositive, {\it i.e.}, $\A_k$ is a symmetric 
$Z$-matrix $(0 \leq k \leq m)$). \vspace{-2mm}
\item{(ii) } The graph $G(L,\FC)$ is a forest and $\sigma_{ij} \in \{ -1,1\}$ 
for every $(i,j) \in \FC$.\vspace{-2mm}
\item{(iii) } The graph $G(L,\FC)$ is bipartite and $\sigma_{ij} = 1$ 
for every $(i,j) \in \FC$.\vspace{-2mm}
\end{description}
Then sub-SDP~\eqref{eq:subSDP1}  is exact for every $\bdelta\in\Real^m$ and $\U \in \SymMat^{\ell}(\DC)$. 
The condition (i) originally was proposed in \cite{KIM2003} for SDP exactness. 
\ecoro

The theorem and corollary above establish exactness for every local
sub-SDP in this sign-pattern class. Thus, under the diagonal consistency 
assumption~\eqref{eq:blockClique}, this class is
parameter-independent: once the sign-pattern 
conditions are satisfied, no further assumption is required on the right-hand-side vector $\bdelta$ or the consistency matrix $\U$.

\subsection{Separable sub-QCQPs characterized by a limited number of constraints}

\label{section:separable}

We now consider a class of separable sub-QCQPs whose objective and
constraint functions share a separable structure.  In contrast to the  classes considered in 
Sections~\ref{section:convex} and~\ref{section:signPattern},  
exactness of the corresponding SDP relaxations is not
guaranteed by  structural properties. 
It depends on the number of
constraints and the right-hand-side vector $\bdelta$ 
induced by the global problem.  Exactness of the SDP
relaxation for this class follows from a rank argument  
(\cite[Theorem~2.2]{PATAKI1998}) that depends
on the number of constraints relative to the number of separable blocks. 
This dependence explains why the class must be  treated separately: 
when a separable subproblem appears
as a  clique-wise component  of the global SDP, the relevant right-hand-side vector $\bdelta$ 
is  not chosen independently, but is induced by an optimal
solution of the  global problem.

\begin{figure}[t!]  
\mbox{ \ } \hspace{10mm}
$
\begin{pmatrix} 
1  &     &   &     &       &     &  \\
   &  2   &   &    &       &     &   \\
   &     & 3 & *   &     &       &    \\
   &     & *  & 4  &     &      &    \\
   &     &    &     & 5   &  *  &  *  \\
   &     &    &     &  *   & 6  &   *       \\
   &     &    &     &  *   &  *   & 7    \\
\end{pmatrix}
$

\vspace{-35mm}
\mbox{ \ } \hspace{70mm}
\includegraphics[height=25mm]{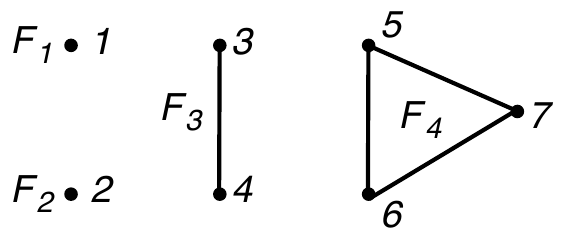}

\vspace{10mm}

\caption{
An example of the aggregate sparsity pattern matrix of a block 
diagonal (or separable) QCQP and its aggregate sparsity pattern graph 
$G(L,\FC)$ with node set $L = \{1,\ldots,7\}$ and edge set 
$
 \FC = \{(3,4),(5,6),(5,7),(6,7)\}. 
$
The maximal cliques are
$F_1=\{1\},F_2=\{2\},F_3=\{3,4\}$ and $F_4=\{5,6,7\}$, which are disjoint.
}
\end{figure} 

To describe the sub-SDP corresponding to a separable sub-QCQP, we impose the
following condition on sub-SDP~\eqref{eq:subSDP1}: \vspace{-2mm} 
\begin{description}
\item{(A) } All the data matrices $\A_k$ $(k=0,\ldots,m)$ 
share a common block-diagonal  structure. 
 We denote the $q$th diagonal block of 
each $\A_k \in \SymMat^{\ell}$ by 
$\B^q_k = \A_k^{F_q} \in \SymMat^{F_q}$ and the corresponding $q$th diagonal block of 
the variable matrix $\Y \in \SymMat^{\ell}$ by $\W^q = \Y^{F_q} \in \SymMat^{F_q}$ $(q=1,\ldots,\hat{q})$,
 where $\{ F_q: q=1,\ldots,\hat{q}\}$ is 
a partition of $C = \{1,\ldots,\ell\}$ such that $i < j$ if $i \in F_q$, $j \in F_r$ and $q < r$. 
 We simply write $\A_k = {\rm diag}(\B^1_k,\ldots,\B^{\hat{q}}_k)$ 
 and $\Y = {\rm diag}(\W^1,\ldots,\W^{\hat{q}})$.
\vspace{-2mm}
\end{description}
In this case, sub-SDP~\eqref{eq:subSDP1} possesses  a block-diagonal sparsity structure. 
Its aggregate sparsity pattern graph itself is  chordal  and  
consists of $\hat{q}$ disjoint cliques $F_q$ $(q=1,\ldots,\hat{q})$. 
See Figure 2, where  
$\ell = 7$, $\hat{q} = 4$, $F_1 = \{1\}, F_2 = \{2\}, F_3 = \{3,4\}$ and $F_4 = \{5,6,7\}$. 

Let $\DC_q = \DC\cap(F_q\times F_q),\ (q=1,\ldots,\hat{q})$, where 
the consistency set $\DC \subseteq \{(i,i) : 1 \leq i \leq \ell\}$ 
is given for sub-SDP~\eqref{eq:subSDP1}. 
Then we can convert 
sub-SDP~\eqref{eq:subSDP1}   
further into a clique-wise SDP as follows. 
\begin{eqnarray}
\eta(\bdelta,\U)
=
\inf\left\{
\displaystyle \sum_{q=1}^{\hat{q}} \inprod{\B^q_0}{\W^q} :
\begin{array}{l}
\W^q \in \SymMat^{F_q}_+
\ (q=1,\ldots,\hat{q}), \\
\displaystyle \sum_{q=1}^{\hat{q}} \inprod{\B^q_k}{\W^q}
\trianglelefteq_k \delta_k \ (k=1,\ldots,m), \\
{[\W^q]_{ii}}=U_{ii}
\ ((i,i)\in \DC_q,\ q=1,\ldots,\hat{q})
\end{array}
\right\}.
\label{eq:subSDP4}
\end{eqnarray}
Under Condition~(A), Lemma~\ref{lemma:chordalGraph} shows that
sub-SDP~\eqref{eq:subSDP1} is equivalent to SDP~\eqref{eq:subSDP4}.
SDP~\eqref{eq:subSDP4} is interpreted as a separable problem because the
matrix variables $\W^q \in \SymMat^{F_q}$ $(q=1,\ldots,\hat q)$ do not
overlap, {\it i.e.}, $F_q\cap F_r=\emptyset$ $(q\ne r)$.

\theo \label{theorem:constLimited}
Assume that Condition~(A) is satisfied, so that
sub-SDP~\eqref{eq:subSDP1} is equivalent to SDP~\eqref{eq:subSDP4}.
Let $\bdelta \in \Real^m$ and $\U \in \SymMat^{\ell}(\DC)$.
Define
\[
\mu =
\mbox{the number of elements in }
\left\{(i,i)\in \DC : U_{ii} > 0\right\}.
\]
Assume that
\vspace{-2mm}
\begin{description}
\item{(B)}
For every optimal solution $(\W^1,\ldots,\W^{\hat q})$ of
SDP~\eqref{eq:subSDP4}, at least $m+\mu-1$ members of the following
collection are nonzero:
\[
\{\W^q : q=1,\ldots,\hat q\}
\ \cup\
\left\{
\delta_k - \sum_{q=1}^{\hat q}\inprod{\B^q_k}{\W^q}
: k=1,\ldots,m
\right\}.
\]
\vspace{-2mm}
\end{description}
Then SDP~\eqref{eq:subSDP1} is exact.\vspace{-2mm}
\etheo
\proof{
By the definition of exactness, it suffices to show
that if sub-SDP~\eqref{eq:subSDP1} has an optimal solution, then it has
a rank-at-most-one optimal solution.
Assume that sub-SDP~\eqref{eq:subSDP1} has an optimal solution.  Under
Condition~(A), sub-SDP~\eqref{eq:subSDP1} is equivalent to
SDP~\eqref{eq:subSDP4}; hence SDP~\eqref{eq:subSDP4} also has an
optimal solution.  Introducing auxiliary variables
$s_k\in\Real$ $(k=1,\ldots,m)$, we rewrite SDP~\eqref{eq:subSDP4} as
\begin{eqnarray}
\eta(\bdelta,\U)
=
\inf\left\{
\displaystyle \sum_{q=1}^{\hat q}\inprod{\B^q_0}{\W^q} :
\begin{array}{l}
\W^q\in\SymMat^{F_q}_+ \quad (q=1,\ldots,\hat q), \\
0\trianglelefteq_k s_k \quad (k=1,\ldots,m), \\
\displaystyle
\sum_{q=1}^{\hat q}\inprod{\B^q_k}{\W^q}+s_k
=
\delta_k \quad (k=1,\ldots,m), \\
{[\W^q]_{ii}}=U_{ii}
\quad ((i,i)\in\DC_q,\ q=1,\ldots,\hat q)
\end{array}
\right\}.
\label{eq:subSDP5}
\end{eqnarray}
Here
$ 
s_k
=
\delta_k-\sum_{q=1}^{\hat q}\inprod{\B^q_k}{\W^q}
\ (k=1,\ldots,m).
$ 

If $U_{ii}=0$ for some equality constraint
$[\W^q]_{ii}=U_{ii}$, then positive semidefiniteness of
$\W^q$ implies that the $i$th row and column of $\W^q$ are zero.
Thus we may eliminate the $i$th row and column from $\W^q$ and
$\B^q_k$ $(k=0,\ldots,m)$, and also eliminate $(i,i)$ from $\DC$.
After carrying out this reduction for all such diagonal constraints,
the remaining SDP~\eqref{eq:subSDP5} has $m$ linear equations  and $\mu$ positive diagonal equality 
constraints.

Since the reduced SDP has an optimal solution, we can apply Theorem~2.2 of
\cite{PATAKI1998}, which provides a rank bound for feasible solutions of an SDP, 
for the existence of an optimal solution
$(\widetilde{\W}^1,\ldots,\widetilde{\W}^{\hat q},
\bar s_1,\ldots,\bar s_m)$ satisfying
\begin{equation}
\sum_{q\in Q}
\frac{\rank(\widetilde{\W}^q)(\rank(\widetilde{\W}^q)+1)}{2}
+
\sum_{k\in K}1
\leq
m+\mu,
\label{eq:rankCond}
\end{equation}
where
$ 
Q=\{q:\widetilde{\W}^q\ne\O\} \ \mbox{and} \ 
K=\{k:\bar s_k\ne0\}
$.
Suppose, to the contrary, that
$\rank(\widetilde{\W}^r)\geq2$ for some $r\in Q$.  Then
\[
\frac{\rank(\widetilde{\W}^r)
(\rank(\widetilde{\W}^r)+1)}{2}
\geq 3.
\]
By Condition~(B), at least $m+\mu-1$ elements among the matrices
$\widetilde{\W}^q$ and the residuals
$ 
\bar s_k = \delta_k-\sum_{q=1}^{\hat q}\inprod{\B^q_k}{\widetilde{\W}^q}
 \  (k=1,\ldots,m)
$ 
are nonzero.  Hence
$ 
|Q\setminus\{r\}|+|K|\geq m+\mu-2.
$ 
Therefore,
\begin{eqnarray*}
m+\mu+1
&=&
3+(m+\mu-2) \\
&\leq&
\frac{\rank(\widetilde{\W}^r)
(\rank(\widetilde{\W}^r)+1)}{2}
+
\sum_{q\in Q\setminus\{r\}}
\frac{\rank(\widetilde{\W}^q)(\rank(\widetilde{\W}^q)+1)}{2}
+
\sum_{k\in K}1,
\end{eqnarray*}
which contradicts~\eqref{eq:rankCond}.  Thus
$\rank(\widetilde{\W}^q)\leq1$ $(q=1,\ldots,\hat q)$.

By Lemma~\ref{lemma:chordalGraph}(ii), there exists a
rank-at-most-one completion
$\widetilde{\Y}\in\SymMat^\ell_+$ satisfying
$\widetilde{\Y}^{F_q}=\widetilde{\W}^q$ 
$(q=1,\ldots,\hat q)$.  Since SDP~\eqref{eq:subSDP4} is equivalent to
sub-SDP~\eqref{eq:subSDP1} under Condition~(A), this
$\widetilde{\Y}$ is a rank-at-most-one optimal solution of
sub-SDP~\eqref{eq:subSDP1}.  Hence SDP~\eqref{eq:subSDP1} is exact.
\qed
}

\medskip

It should be noted that Condition (B) involves  
$\U \in \SymMat^{\ell}(\DC)$, which determines $\mu$, and 
$\bdelta = (\delta_1,\ldots,\delta_m) \in \Real^m$. 
Thus, exactness of 
sub-SDP~\eqref{eq:subSDP1}
depends  not only on the number $m$ of constraints, but also on the right-hand-side vector $\bdelta \in \Real^m$  and the  consistency matrix $\U \in \SymMat^{\ell}(\DC)$. 
In Theorem~\ref{theorem:main1}, these are not free parameters: they
are induced by an optimal solution 
$(\widetilde{\Y}^1,\ldots,\widetilde{\Y}^{\hat{p}},\widetilde{\U})$
of the global  SDP~\eqref{eq:SDP4},
with $\bdelta$ corresponding to 
one of the $\tilde{\bdelta}^p$ $(p=1,\ldots,\hat{p})$, and $\U$ corresponding to  $\widetilde{\U}$.
Consequently, when SDP~\eqref{eq:subSDP4} appears as a local sub-SDP within the global SDP, 
its exactness may depend on the  behavior of the other local sub-SDPs.
This dependence distinguishes the present separable class from the parameter-independent classes 
based on convexity or sign-pattern conditions.

\rema
Theorem~\ref{theorem:constLimited} can be compared with the result in \cite{LUO2010}, 
where the following assumptions were imposed:
\begin{eqnarray*}
& & m \leq {\hat{q}}+1, \ \DC_q = \emptyset \ (q=1,\ldots,\hat{q}) \ (\mbox{hence $\mu = 0$}), \\
& & \W^q \not= \O \ (1 \leq q \leq \hat{q}) \ \mbox{for every optimal solution } (\W^1,\ldots,\W^{\hat{q}}).
\end{eqnarray*}
Under these assumptions, \eqref{eq:rankCond} implies 
\begin{eqnarray*}
m-1 \leq {\hat{q}} \leq \sum_{q=1}^{\hat{q}}
\frac{\mbox{rank}(\widetilde{\W}^q)(\mbox{rank}(\widetilde{\W}^q)+1)}{2} 
\leq m,
\end{eqnarray*}
and hence either $\hat q = m-1$ or $\hat q = m$, 
that is, either $m=\hat q+1$ or $m=\hat q$. 
In contrast, Condition~(B) of Theorem~\ref{theorem:constLimited} 
considerably relaxes this restriction.  
\erema

Even in case $\DC_q=\emptyset$ $(q =1,\ldots,\hat{q})$, whether Condition~(B) holds depends 
on the right-hand-side vector $\bdelta$  
and the relations `$\trianglelefteq_k$' $(1 \le k \le m)$. 
When $m  + |\DC| \leq 2$, the conclusion of Theorem~\ref{theorem:constLimited} holds 
without any additional assumptions on $\bdelta$, 
the relations `$\trianglelefteq_k$', and $\B^q_k$ 
$(q=1,\ldots, {\hat{q}}, \ k=0,\ldots,m)$, and the number $\hat{q}$ of the separable blocks 
as shown below. (This result is known; see \cite{ARIMA2023,POLYAK1998,YE2003}.) 
\coro \label{corollary:m=2} 
Assume that $m + |\DC| \leq 2$.\vspace{-2mm} 
\begin{description}
\item{(i) } SDP~\eqref{eq:subSDP4} is exact. 
\vspace{-2mm} 
\item{(ii) } SDP~\eqref{eq:subSDP1} is exact. 
\vspace{-2mm} 
\end{description}
\ecoro
\proof{ 
Since (ii) follows from (i) by taking the single block
$F_1=C$, we only prove (i).  If Condition~(B) of
Theorem~\ref{theorem:constLimited} holds, then the conclusion follows
immediately.
Suppose that Condition~(B) fails.  Since
$\mu\leq |\DC|$ and $m+|\DC|\leq 2$, we have
$m+\mu\leq 2$.  If $m+\mu\leq 1$, then Condition~(B) cannot fail.
Hence the failure of Condition~(B) implies $m+\mu=2$, and there exists
an optimal solution
$(\widetilde{\W}^1,\ldots,\widetilde{\W}^{\hat q})$ of
SDP~\eqref{eq:subSDP4} for which at most $m+\mu-2=0$ elements among
the matrices $\widetilde{\W}^q$ and the residuals are nonzero.
Therefore $\widetilde{\W}^q=\O$ $(q=1,\ldots,\hat q)$, and hence
$\rank(\widetilde{\W}^q)=0\leq 1$ $(q=1,\ldots,\hat q)$.
This proves the result.
\qed
}

\subsection{Preserving the exactness of sub-SDP~\eqref{eq:subSDP1} by 
adding dependent inequality constraints} 

The following theorem establishes a preservation result for a family of
sub-SDPs with varying right-hand-side vectors.  
Assume that sub-SDP~\eqref{eq:subSDP1}
has a rank-at-most-one optimal solution for every right-hand-side vector whenever
 its optimal value is finite and attained.  
 Then inequality constraints whose coefficient matrices lie in the conic hull of 
the coefficient matrices of the existing constraints 
may be added without destroying exactness. 
Thus, the resulting modified sub-SDP is exact.

\theo \label{theorem:extension}
Assume that $\trianglelefteq_k=$ `$\le$' $(k=1,\ldots,m)$ and
$\U\in\SymMat^\ell(\DC)$.  
Suppose that, for every $\bdelta'\in\Real^m$, if
$-\infty < \eta(\bdelta',\U) < \infty$, then
 SDP~\eqref{eq:subSDP1} with $\bdelta=\bdelta'$ 
has a rank-at-most-one optimal solution.
Let $\bdelta\in\Real^m$ and $\delta_j\in\Real$
$(j=m+1,\ldots,m')$ be fixed, where $m < m'$.   Assume that
\begin{eqnarray*}
& & \A_j\in
\mbox{\rm cone}\{\A_k:k=1,\ldots,m\} \equiv
\left\{
\sum_{k=1}^m \alpha_k\A_k :
\alpha_k \geq 0 \ (k=1,\ldots,m)
\right\} \\
& &  \hspace{90mm} (j=m+1,\ldots,m').
\end{eqnarray*}
Then  the modified SDP obtained from SDP~\eqref{eq:subSDP1}
by adding the constraints
\[
\inprod{\A_j}{\Y}\leq \delta_j
\ (j=m+1,\ldots,m')
\]
is exact. 
\etheo
\proof{
Let $\widetilde{\Y}\in\SymMat^\ell_+$ be an optimal solution of the
modified SDP.  Define
\[
\delta'_k=\inprod{\A_k}{\widetilde{\Y}} \ (k=1,\ldots,m), \ \bdelta'=(\delta'_1,\ldots,\delta'_m).
\]
Then $\delta'_k\leq\delta_k$ $(k=1,\ldots,m)$, and
$\widetilde{\Y}$ is feasible for SDP~\eqref{eq:subSDP1} with
right-hand-side vector $\bdelta'$.
For each $j=m+1,\ldots,m'$, take coefficients
$\alpha_{jk}\geq0$ $(k=1,\ldots,m)$ such that
$ 
\A_j=\sum_{k=1}^m \alpha_{jk}\A_k .
$ 
We show that all the added constraints are redundant for
SDP~\eqref{eq:subSDP1} with right-hand-side vector $\bdelta'$.
Indeed, if $\Y\in\SymMat^\ell_+$ satisfies
$ 
\inprod{\A_k}{\Y}\leq \delta'_k
\quad (k=1,\ldots,m),
$ 
then, for every $j=m+1,\ldots,m'$,
\[
\inprod{\A_j}{\Y}
=
\sum_{k=1}^m\alpha_{jk}\inprod{\A_k}{\Y}
\leq
\sum_{k=1}^m\alpha_{jk}\delta'_k
=
\sum_{k=1}^m\alpha_{jk}\inprod{\A_k}{\widetilde{\Y}}
=
\inprod{\A_j}{\widetilde{\Y}}
\leq
\delta_j.
\]
Thus, after replacing $\bdelta$ by
$\bdelta'$,  the added constraints are redundant.
Consequently, SDP~\eqref{eq:subSDP1} with
right-hand-side vector $\bdelta'$ and the modified SDP have the same optimal value. 
Indeed, every feasible solution of the former  is
feasible for the modified SDP by the redundancy shown above, while
$\widetilde{\Y}$ is feasible for the former and attains the
optimal value of the modified SDP.
Therefore
$ 
-\infty < \eta(\bdelta',\U) < \infty .
$ 
By the assumption, SDP~\eqref{eq:subSDP1} with right-hand-side vector
$\bdelta'$ has a rank-at-most-one optimal solution.  This solution is
also feasible for the modified SDP and has the same optimal value.
Hence it is a rank-at-most-one optimal solution of the modified SDP.
\qed
}

The proof above shows that, although the additional constraints may
shrink the feasible region, they become redundant after replacing the
right-hand-side vector by the values attained at an optimal solution of
the modified SDP.  This preservation result is particularly useful when the resulting
SDP appears as a local sub-SDP in the clique-wise relaxation
\eqref{eq:SDP4}. In that setting,  
 the relevant local right-hand side vector is induced by an optimal
solution of the global SDP. 
Combining Theorem~\ref{theorem:extension} with Corollary~\ref{corollary:m=2}(ii), 
we obtain the following extension.

\coro \label{corollary:dependent}
Assume that $m\in\{1,2\}$ and $m+|\DC|\leq 2$ in
SDP~\eqref{eq:subSDP1}.  Let
\[
\A_j\in \mbox{\rm cone}\{\A_k:k=1,\ldots,m\},
\
\delta_j\in\Real
\ (j=m+1,\ldots,m'), 
\]
where $m <m'$. 
Then the modified SDP obtained from SDP~\eqref{eq:subSDP1}
by adding the constraints
$ 
\inprod{\A_j}{\Y}\leq \delta_j
\ (j=m+1,\ldots,m')
$ 
is exact. 
\ecoro

We call the inequalities
$ 
\inprod{\A_j}{\Y}\leq \delta_j
\ (j=m+1,\ldots,m')
$ 
{\em dependent inequalities} if their coefficient matrices satisfy
$ 
\A_j\in \mbox{\rm cone}\{\A_k:k=1,\ldots,m\}
\ (j=m+1,\ldots,m').
$ 
The dependence occurs only at the level of the coefficient matrices; 
the right-hand-side constants $\delta_j$ $(j=m+1,\ldots,m')$ are arbitrary real numbers.
Here $\A_1,\ldots,\A_m$ serve as the {\em base coefficient matrices}, and the
remaining matrices $\A_j$ $(j=m+1,\ldots,m')$ are dependent on them 
in the sense that they lie in the conic hull of 
the base coefficient matrices.
For notational convenience, we call all inequalities
$\inprod{\A_j}{\Y}\leq\delta_j$ $(j=1,\ldots,m')$  {\em a system of
dependent inequalities} generated by  $\A_1,\ldots,\A_m$.

\section{Examples}

\label{section:example}

In this section, we present three examples that illustrate how the
exactness of the SDP relaxation~\eqref{eq:SDP0} of QCQP~\eqref{eq:QCQP0}
can be derived from local sub-SDPs by combining
Theorem~\ref{theorem:main1} with the local exactness results in
Section~\ref{section:subQCQPs}.  These local results include 
sub-QCQPs characterized by convexity, sign-pattern conditions and a system of 
dependent inequalities.  The final example shows how
heterogeneous classes of sub-QCQPs including separable sub-QCQPs characterized 
by a limited number of constraints can be incorporated 
within the local-to-global exactness framework while preserving exactness of the SDP
relaxation.

\examp \label{example:separable} 
This example illustrates the simplest case in which the clique-wise
reformulation becomes completely separable and no consistency constraints
are required.  We assume that $\trianglelefteq_k=$ `$\leq$'
$(k=1,\ldots,m)$ and that each $\A_k$ in QCQP~\eqref{eq:QCQP0} and its
SDP relaxation~\eqref{eq:SDP0} is of the following block diagonal form:
\begin{eqnarray*}
\A_k =
\begin{pmatrix}
\A^1_k & \O & \O \\
\O & \A^2_k & \O \\
\O & \O & \A^3_k
\end{pmatrix}
\in \SymMat^n, \quad
\A^p_k \in \SymMat^{C_p}
\quad (p=1,2,3,\ k=0,\ldots,m),
\end{eqnarray*}
where
$C_1=\{1,\ldots,n_1\}$,
$C_2=\{n_1+1,\ldots,n_2\}$,
$C_3=\{n_2+1,\ldots,n_3\}$, and $n_3=n$.
Since $C_p\cap C_q=\emptyset$ $(p\ne q)$, we have
$\DC_p=\emptyset$ $(p=1,2,3)$ and $\overline{\DC}=\emptyset$.
The aggregate sparsity pattern graph consists of three disconnected
cliques $C_p$ $(p=1,2,3)$.
Hence the clique-wise formulation~\eqref{eq:SDP4} of SDP~\eqref{eq:SDP0}
becomes separable.

If we assign to each clique $C_p$ a sub-QCQP characterized by convexity
(Theorem~\ref{theorem:convex}), sign-pattern conditions
(Theorem~\ref{theorem:signDefinite} and
Corollary~\ref{corollary:signDefinite}), 
or a system of dependent inequalities (Corollary~\ref{corollary:dependent}), 
then
SDP~\eqref{eq:SDP4}, and hence SDP~\eqref{eq:SDP0}, is exact.
For example, we can take
$\A^1_k\in\SymMat^{C_1}_+$ $(k=0,\ldots,m)$ for a convex sub-QCQP on
$C_1$, and $\A^2_k\in\SymMat^{C_2}$ with all off-diagonal elements
nonpositive $(k=0,\ldots,m)$ for a sub-QCQP characterized by
sign-pattern conditions (Corollary~\ref{corollary:signDefinite}). 
For the sub-QCQP on $C_3$ with dependent inequalities, assume $m\geq2$ and
take two base coefficient matrices
$ 
\O\ne \A^3_1\in\SymMat^{C_3}, \
\O\ne \A^3_2\in\SymMat^{C_3}.
$ 
Assume that the remaining coefficient matrices satisfy
$ 
\A^3_k\in\mbox{\rm cone}\{\A^3_1,\A^3_2\}
\ (k=3,\ldots,m),
$ 
while $\A^3_0\in\SymMat^{C_3}$ is arbitrary.  Since $\DC_3=\emptyset$,
Corollary~\ref{corollary:dependent} applies to this local sub-SDP.

The discussion after SDP~\eqref{eq:SDP0} also applies to the present
example: linear terms can be incorporated by including, or adding, a
normalization constraint $X_{ii}=1$ for some $i \in \{1,\ldots,n_2\}$.  The same observation applies to
Examples~\ref{example:52} and~\ref{example:53}.  If this normalization is
added as an extra constraint, however, the number of equality constraints 
increases.
This must be taken into account when applying results whose assumptions
depend on the number of equality and/or inequality constraints, such as Theorem~\ref{theorem:convex},
Theorem~\ref{theorem:constLimited}
and Corollary~\ref{corollary:dependent}.
\eexamp

\examp \label{example:52}
This example illustrates how exactness can be preserved under a
clique-wise coupling of sub-QCQPs characterized by convexity and those
characterized by sign-pattern conditions, with a nonempty consistency set 
$\overline{\DC}$. 
We consider QCQP~\eqref{eq:QCQP0} with $n=5$, $m \geq 1$, 
$\trianglelefteq_k=$ `$\le$'  $(k=1,\ldots,m)$ and $\A_k$ $(k=0,\ldots,m)$ whose 
aggregate sparsity pattern matrix and graph $G(N,\EC^0)$ are given by the following.

\bigskip

\hspace{20mm}
\begin{minipage}{50mm}
$ 
\begin{pmatrix} \oplus  & *                 &                 &               &  \\
                         *          & *                 & \ominus   & \ominus  &      \\
                                    & \ominus      & *               & \ominus  &  \\                         
                                    & \ominus      & \ominus   & *             &  *     \\                         
                                   &                     &                &*               & \oplus                               
                        
\end{pmatrix} 
$ 
\end{minipage}
\begin{minipage}{80mm}
\includegraphics[height=20mm]{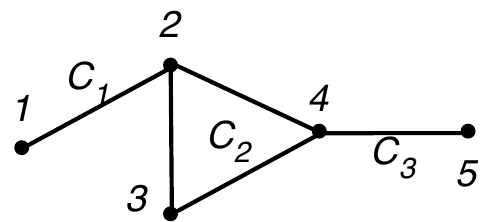}
\end{minipage}

\bigskip

\noindent
Here $\oplus$ denotes a nonnegative real number, $\ominus$ a nonpositive real number and 
$*$ an arbitrary real number. The graph $G(N,\EC^0)$ itself is a block-clique graph with
$\hat{p} = 3$ maximal cliques $C_1= \{1,2\}, C_2 = \{2,3,4\}$ and $C_3 = \{4,5\}$. 
The corresponding consistency sets $\DC_1 = \{(2,2)\}$,  $\DC_2 = \{(2,2),(4,4)\}$ 
and $\DC_3 = \{(4,4)\}$ are diagonal. 
Define
\begin{eqnarray*}
& & \A^p_k = \A^{C_p}_k \in \SymMat^{C_p} \ (k=0,1,\ldots,m,p = 1,3), \\ 
& & \A^2_k = \begin{pmatrix}0                 & [\A_k]_{23} & [\A_k]_{24}  \\  
                                              [\A_k]_{32} & [\A_k]_{33} & [\A_k]_{34}  \\
                                              [\A_k]_{42} & [\A_k]_{43} & 0                                                
                                              \end{pmatrix} \in \SymMat^{C_2} \ (k=0,1,\ldots,m), 
\end{eqnarray*}
so that 
$\A^p_k$ $(k=0,\ldots,m,p=1,2,3)$ satisfy~\eqref{eq:cliqueDecomposition0}. 

\bigskip

Applying the clique-wise reformulation discussed in
Section~\ref{section:clique-wiseReformulation} to the SDP
relaxation~\eqref{eq:SDP0} of QCQP~\eqref{eq:QCQP0}, we obtain
SDP~\eqref{eq:SDP4}.  
The sub-QCQPs on $C_1$ and $C_3$ fall into  the
convex class in Theorem~\ref{theorem:convex}, case (I).  The sub-QCQP on
$C_2$ falls into the sign-pattern class of 
Corollary~\ref{corollary:signDefinite}, case (i), because all its
off-diagonal coefficients are nonpositive.
Therefore, for the parameters $\bdelta^p\in\Real^m$
$(p=1,2,3)$ and the consistency matrix
$\U \in \SymMat^n(\overline{\DC})$ induced by any optimal solution of
SDP~\eqref{eq:SDP4}, all local sub-SDPs are exact.  It follows from
Theorem~\ref{theorem:main1}(iii) that SDP~\eqref{eq:SDP4}, and consequently
SDP~\eqref{eq:SDP0}, is exact.

This example illustrates that exactness can be preserved even when
different local exactness mechanisms are assembled through diagonal 
consistency constraints.
\eexamp

\examp \label{example:53} 
This example illustrates how the local-to-global exactness framework 
can incorporate
three heterogeneous classes of QCQPs with exact SDP relaxations,
including a separable sub-QCQP characterized by a limited
number of constraints discussed in
Theorem~\ref{theorem:constLimited}. 
The main point of this example is the treatment of the separable 
sub-QCQP.  
In contrast to sub-QCQPs characterized by convexity or  sign-pattern conditions, 
the exactness of the associated sub-SDP~\eqref{eq:subSDP4}
 is parameter-dependent;
 it depends on the parameters $(\bdelta,\U)$
induced by the global SDP~\eqref{eq:SDP4}.  
We show that suitable conditions on the other
sub-SDPs can force the sub-SDP~\eqref{eq:subSDP4} to satisfy Condition~(B) of 
Theorem~\ref{theorem:constLimited}.

We construct an instance of QCQP~\eqref{eq:QCQP1} by combining five sub-QCQPs on 
$C_1=\{1,\ldots,6\}$, $C_2=\{5,7\},C_3=\{6,7,8\},C_4=\{7,9\}$ and 
$C_5=\{6,10\}$. See Figure~3.
The node set $C_1=\{1,\ldots,6\}$ involves four cliques $F_1,F_2,F_3,F_4$, 
while $C_2,C_3,C_4,C_5$ are cliques.   
We assign a separable QCQP 
to  $C_1$, QCQPs satisfying sign-pattern conditions to  $C_2$ and $C_3$, convex QCQPs 
to  $C_4$ and $C_5$.  Each sub-QCQP has 
$m=5$ inequality constraints, including some redundant constraints. Table 1 provides the  
details of those sub-QCQPs. 

\begin{figure}[t!]  
\begin{center}
\includegraphics[height=35mm]{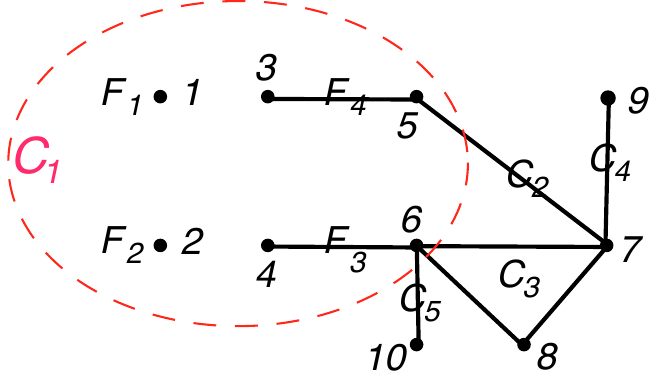}
\caption{
The aggregate sparsity pattern graph $G(N,\EC^0)$ of QCQP~\eqref{eq:QCQP0} 
of Example~\ref{example:53}.
Here $F_q$ $(q=1,\ldots,4)$ and $C_p$ $(p=2,\ldots,5)$ are cliques of $G(N,\EC^0)$, 
whereas $C_1$ is a clique of a chordal extension of the graph. 
Separable QCQPs  with a  limited number of constraints are assigned to 
$C_1 = F_1\cup F_2\cup F_3\cup F_4$, QCQPs characterized by 
sign-pattern conditions to $C_2$ and $C_3$, and convex QCQPs to $C_4$ and $C_5$. 
}
\end{center}
\end{figure} 

\begin{table}[h!] 
\begin{center}
\scriptsize{
\begin{tabular}{|c|c|c|c|c|c||c||c|c|c|c|c|c|} 
\hline
&           &  \multicolumn{4}{|c||}{ $p=1: C_1 = \{1,2,3,4,5,6\}$}                              & $p=2,3$       & $p=4,5$& & \\
&           & \multicolumn{4}{|c||}{Separable} & Sign Pat.
& Convex & & \\
&           & \multicolumn{4}{|c||}{Sect.4.3} & Sect.4.2
& Sect.4.1 & & \\
\hline
&            & $q=1$ & $q=2$  &  $q=3$ & $q=4$ & & & & \\ 
&            & $F_1=\{1\}$ & $F_2=\{2\}$  &  $F_3=\{3,5\}$ & $F_4=\{4,6\}$ &$C_2=\{5,7\}$ &$C_4=\{7,9\}$ & & \\ 
&            &  &   &   &  & $C_3=\{6,7,8\}$ & $C_5=\{6,10\}$ & & \\ 
\hline
&$\DC$ &  &   & $(5,5)\in\DC_1$  & $(6,6)\in\DC_1$ & $(5,5),(7,7)\in\DC_2$ & $(7,7)\in\DC_4$ & &  \\ 
& &  &   &               &             & $(6,6),(7,7)\in\DC_3$ & $(6,6)\in\DC_5$ & &  \\ 
\hline
&  $k$    &  $\B_k^1$ &  $\B_k^2$ & $\B_k^3$ & $\B_k^4 $ & $\A^2_k$ and  $\A^3_k$ & $\A^4_k$ and  $\A^5_k$ 
  &$\trianglelefteq_k$ & $b_k$ \\ 
\hline 
obj. & $0$    & $\forall$sym  & $\forall$sym   & $\forall$sym  & $\forall$sym  & off-diag$\ominus$ & convex & &  \\ 
\hline 
&  $1$    & $\forall$sym  &  $\oplus$ & $\oplus$ & $\oplus$ & $\oplus$ \& off-diag$\ominus$ & $\oplus$ &$\le$ & $-$ \\ 
&  $2$    & $\oplus$  & $\forall$sym  & $\oplus$ & $\oplus$ & $\oplus$ \& off-diag$\ominus$ & $\oplus$ &$\leq$ & $-$ \\ 
const. &  $3$    & $\forall$sym   & $\forall$sym  & $\forall$sym & $\forall$sym  &  off-diag$\ominus$  & convex &$\leq$ &$\forall$ \\ 
&  $4$    &  $\O$ &  $\O$ & $\O$ & $\O$ &  off-diag$\ominus$ & convex &$\leq$ & $\forall$\\ 
&  $5$    &  $\O$ &  $\O$ & $\O$ & $\O$ &  off-diag$\ominus$ & convex &$\leq$ & $\forall$\\ 
\hline 
\hline
\end{tabular}	
}
\end{center}
\caption{
Summary of data matrices and conditions imposed on the five sub-SDPs on $C_1,\ldots,C_5$. 
$\oplus$: a positive semidefinite symmetric matrix. 
off-diag$\ominus$: a symmetric matrix with all off-diagonal elements nonpositive.
convex: a symmetric matrix $\A$ such that $\inprod{\A}{\y\y^T}$ is convex in $\y$ 
with $y_i$ fixed for some $i$  (see case (I) in Section 4.1). 
$\forall$sym: an arbitrary symmetric matrix. $\O$: a zero matrix. 
}
\end{table}

For a separable sub-QCQP on $C_1$, 
we assume that $C_1$ consists of $4$ disjoint cliques $F_1 = \{1\}, F_2 = \{2\}, 
F_3=\{3,5\}, F_4=\{4,6\}$ as shown in Figure 3. The sub-QCQP on $C_1$  
contains $5$ inequality constraints, where the last two constraints, 
the $4$th and $5$th constraints are redundant. 
These two redundant constraints are included solely to match 
the number $m=5$ of inequality constraints of the QCQP~\eqref{eq:QCQP1},
which also embeds the  sub-QCQPs on $C_2,\ldots,C_5$. 
Hence, when applying
Theorem~\ref{theorem:constLimited} to the separable sub-SDP on $C_1$, we
first remove them and use the reduced formulation
with the three effective inequalities with $k=1,2,3$.
For Condition~(B), we impose
\begin{eqnarray}
& & \A^{p}_k \in \SymMat^{C_p}_+, \ b_k < 0
\ (p=2,3,4,5,\ k=1,2), \label{eq:Cp}\\
& & \B^{q}_1 \in \SymMat^{F_q}_+ \ (q=2,3,4), \
\B^{q}_2 \in \SymMat^{F_q}_+ \ (q=1,3,4).
\label{eq:Fq}
\end{eqnarray}
The condition~\eqref{eq:Cp} ensures that, for every optimal solution $(\Y^1,\ldots,\Y^5,\U)$ of SDP~\eqref{eq:SDP4}, 
\begin{eqnarray*}
\delta^1_k = \inprod{\A^{1}_k}{\Y^1} \le b_k - \sum_{p=2}^5 \inprod{\A^{p}_k}{\Y^p} < 0 \ (k=1,2). 
\end{eqnarray*}
Then the condition~\eqref{eq:Fq} ensures that 
$\W^1 \ne \O$ and $\W^2 \ne \O$ for every optimal solution $(\W^1,\ldots,\W^4)$ of the separable SDP~\eqref{eq:subSDP4} for these $\delta^1_1 < 0 $, $\delta^1_2<0 $ 
and any $\delta^1_3 \in \Real$. 

If both $U_{55}$ and $U_{66}$ are positive, then $\mu=2$ in the reduced
separable sub-SDP.  Since the reduced formulation has three effective
inequality constraints, Condition~(B) requires
$3+\mu-1=4$ nonzero elements.  The conditions above imply
$\W^1\ne\O$ and $\W^2\ne\O$, while $U_{55}>0$ and $U_{66}>0$ imply
$\W^3\ne\O$ and $\W^4\ne\O$.  Hence Condition~(B) is satisfied.

If $U_{55}=0$ and/or $U_{66}=0$, the corresponding row and column can be
eliminated by positive semidefiniteness, as described in
Theorem~\ref{theorem:constLimited}.  The number $\mu$ decreases accordingly,
and the required number of nonzero elements in Condition~(B) decreases by
the same amount.  Therefore, the reduced separable sub-SDP on $C_1$ is exact.

For the sub-QCQPs on $C_2$ and  $C_3$,
we impose the sign-pattern conditions.  Specifically,
 all off-diagonal elements of $\A^{p}_k$
$(k=0,\ldots,5,\ p=2,3)$ are assumed to be nonpositive.  By
Corollary~\ref{corollary:signDefinite}(i), the corresponding
sub-SDPs are exact for every $\bdelta\in\Real^m$ and $\U$.
In addition, to ensure  that the separable sub-QCQP on $C_1$ satisfies Condition (B), 
we have assumed that $\A^{p}_k \in \SymMat^{C_p}_+$ $(k=1,2,p=2,3)$.

For the sub-QCQPs on $C_4$ and $C_5$, 
we use the convex QCQPs discussed in Section~\ref{section:convex}. 
Since $\DC_4=\{(i_4,i_4)\}=\{(7,7)\}$ and
$\DC_5=\{(i_5,i_5)\}=\{(6,6)\}$, case (I) of
Theorem~\ref{theorem:convex} applies to both sub-QCQPs.
Thus, we impose the convexity condition that 
the principal submatrix of $\A^p_k$ indexed by 
$C_p\backslash\{i_p\}$ 
is positive semidefinite $(k=0,\ldots,5,p=4,5)$, or equivalently, 
$ 
[\A^4_k]_{99} \geq 0, \ [\A^5_k]_{10,10} \geq 0 \ (k=0,\ldots,5). 
$ 
It follows  that  the sub-SDPs~\eqref{eq:subSDP1} on $C_4$ and  $C_5$ are exact 
for every $\bdelta\in\Real^m$ and $\U$. 
In order for the separable sub-QCQP on $C_1$ 
to satisfy Condition (B), we have assumed  that $\A^{p}_k \in \SymMat^{C_p}_+$ 
$(k=1,2,p=4,5)$. 

Consequently,  each local sub-SDP induced by an optimal solution of
SDP~\eqref{eq:SDP4} is exact.   It follows from 
Theorem~\ref{theorem:main1}(iii) that SDP~\eqref{eq:SDP4}, and therefore
SDP~\eqref{eq:SDP0}, is exact.
This example highlights a feature that does not arise in the purely
convex or sign-pattern cases: the exactness of one local sub-SDP may be
certified using parameter information imposed by the other  local
sub-SDPs through the global constraints.
\eexamp

\section{Concluding remarks}

\label{section:conclusion}

We have developed a local-to-global exactness framework for 
SDP relaxations of sparse QCQPs.  Using a chordal extension of the
aggregate sparsity pattern graph, the SDP relaxation is reformulated in
clique-wise matrix variables associated with the maximal cliques.  This
equivalent reformulation induces local sub-SDPs linked by
consistency constraints on clique overlaps.  
The main result shows that exactness of the original SDP relaxation can
be certified by rank-at-most-one attainability of the local sub-SDPs with
the local right-hand-side vectors and consistency matrix induced by an
optimal solution of the clique-wise formulation.

For the applications developed in this paper, the block-clique assumption 
is imposed  in 
handling the consistency of rank-at-most-one optimal solutions of the local sub-SDPs.
Under this assumption, all consistency constraints on clique
overlaps are diagonal, so local rank-at-most-one optimal solutions can be
combined without imposing additional off-diagonal rank-one consistency
conditions.  This allows different local exactness mechanisms to be used
within a single sparse QCQP.  
In particular, we have identified three mechanisms based, respectively on convexity, sign-pattern conditions, 
and separability with a limited number of constraints.
The examples in Section~\ref{section:example} illustrate how these
heterogeneous local exactness mechanisms  
can be assembled to prove exactness of the
original SDP relaxation.

An important issue for future work is the choice of chordal extension
and clique-wise decomposition of the data matrices.  Different chordal
extensions, and different decompositions over their maximal cliques, may
lead to different local sub-SDPs.  Consequently, they may affect whether
the local exactness certificates developed in this paper can be applied.

It would  also be 
valuable to further enlarge the list of local sub-QCQP
classes with exact SDP relaxations.  The framework is modular: any new
local exactness result that is stable under the induced right-hand side
and consistency matrix can be incorporated into the global
certification scheme.  This suggests that sparse exactness analysis may
serve as a way to combine otherwise separate exactness criteria for
nonconvex QCQPs.

\vspace{0.5cm}

\subsection*{Data availability statement} Data sharing is not applicable to this article as no datasets were generated or 
analyzed during the current study.

\subsection*{Declarations} 

{\bf Competing interests} \ The authors have no competing interests to declare that are relevant to the content of this article.


\end{document}